\documentclass{amsart}
\usepackage{amsmath,latexsym}
\usepackage[mathscr]{eucal}
\usepackage{amssymb}

\newtheorem{theorem}{Theorem}

\newtheorem*{Prop1prime}{Proposition 1$'$}

\newtheorem*{Prop2prime}{Proposition 2$'$}

\newtheorem*{Prop1simple}{Proposition 1 (Special Case)}
\newtheorem{proposition}{Proposition}
\newtheorem{lemma}{Lemma}
{Corollary}
\theoremstyle{definition}

\theoremstyle{remark}

\newtheorem*{rema}{Remark}
\numberwithin{equation}{section}

\let\epsilon\varepsilon

\def\sumprime_#1{\setbox0=\hbox{$\scriptstyle{#1}$}
\setbox2=\hbox{$\displaystyle{\sum}$}
\setbox4=\hbox{${}'
\mathsurround=0pt$}
\dimen0=.5\wd0 \advance\dimen0 by-.5\wd2
\ifdim\dimen0>0pt
\ifdim\dimen0>\wd4 \kern\wd4 \else\kern\dimen0\fi\fi
\mathop{{\sum}'}_{\kern-\wd4 #1}}
\newcommand{\gs}{\mathfrak{S}}
\def\ndiv{\not \hskip .03in \mid}
\begin{document}

\title{Primes in Tuples I}
\author{D. A. Goldston}
\address{Department of Mathematics, San Jose
State University, San Jose, CA 95192, USA}
\email{goldston@math.sjsu.edu}
\author{J. Pintz}
\address{R\'enyi Mathematical Institute of the Hungarian Academy
of Sciences, H-1364 Budapest, P.O.B. 127, Hungary} 
\email{pintz@renyi.hu}
\author{C. Y. Y{\i}ld{\i}r{\i}m}
\address{ Department of Mathematics, Bo\~{g}azi\c{c}i University,
Istanbul 34342 \& Feza G\"{u}rsey 
Enstit\"{u}s\"{u},
\c{C}engelk\"{o}y, Istanbul, P.K. 6, 81220, Turkey}
\email{yalciny@boun.edu.tr}
\thanks{The first author was supported by NSF grant DMS-0300563,
the NSF Focused Research Group grant 0244660, and the 
American Institute of Mathematics;
the second author by OTKA grants No. T38396, T43623, T49693 and 
the Balaton program; 
the third author by T\"UB\.{I}TAK}

\begin{abstract}  We introduce a method for showing that there exist prime numbers which are very close together. The method depends on the level of distribution of primes in arithmetic progressions. Assuming the Elliott-Halberstam conjecture, we prove that there are infinitely often primes differing by 16 or less. Even a much weaker conjecture implies that there are infinitely often primes a bounded distance apart. Unconditionally, we prove that there exist consecutive primes which are closer than any arbitrarily small multiple of the average spacing, that is,
\[ \liminf_{n\to \infty} \frac{p_{n+1}-p_n}{\log p_n} =0 .\] 
This last result will be considerably improved  in a later paper.
 \end{abstract}

\maketitle

\section{Introduction}
One of the most important unsolved problems in number theory is to establish the existence of infinitely many prime tuples. Not only is this problem believed to be difficult, but it  has also earned the reputation among most mathematicians in the field as hopeless in the sense that there is no known unconditional approach for tackling the problem. The purpose of this paper, the first in a series, is to provide what we believe is a method which could lead to a partial solution for this problem. At present, our results on primes in tuples are conditional on information about the distribution of primes in arithmetic progressions. However, the information needed to prove that there are infinitely often two primes in a given $k$-tuple for sufficiently large $k$ does not seem to be too far beyond the currently known results.   Moreover, we can gain enough in the argument by averaging over many tuples to obtain unconditional results concerning small gaps between primes which go far beyond anything that has been proved before. Thus, we are able to prove the existence of very small gaps between primes which, however, go slowly to infinity with the size of the primes. 

The information on primes we utilize in our method is often referred
to as the level of distribution of primes in
arithmetic progressions. Let
\begin{equation}
\theta(n) = \left\{
\begin{array}{ll}
  \log n &\text{if } n \text{ prime},\\
  0 &\text{otherwise},
\end{array}
\right. \label{1.1}
\end{equation}
and consider the counting function
\begin{equation}
\theta(N;q,a) = \sum_{\substack{n\le N \\ n\equiv a (\text{mod}\,
q)}}\theta(n). \label{1.2}
\end{equation}
The Bombieri-Vinogradov theorem states that for any  $A>0$ there is a $B=B(A)$ such that, for $Q=N^{\frac{1}{2}}(\log N)^{-B}$,
\begin{equation}
\sum_{q\le Q}  \max_{\substack{a\\  (a,q)=1}}  \left| \theta(N;q,a)-   \frac{N}{\phi(q)}  \right| \ll \frac{N}{(\log N)^A}. \label{1.3}
\end{equation}
  We say that the primes have \emph{level of distribution} $\vartheta$ if \eqref{1.3} holds for any $A>0$ and any $\epsilon >0$ with
\begin{equation} Q=N^{\vartheta -\epsilon}.  \label{1.4} \end{equation}  
Elliott and Halberstam \cite{EH} conjectured that the primes have level of distribution 1. According to the Bombieri-Vinogradov theorem, the primes are known to have level of distribution $1/2$.
 
Let $n$ be a natural number  and consider the $k$-tuple 
\begin{equation}  (n+h_1,n+h_2, \ldots , n+h_k),   \label{1.5} \end{equation}
where $\mathcal{H}=\{h_1,h_2,\ldots , h_k\}$ is a set composed of   distinct  non-negative integers.  If every component of the tuple is a prime we call this a \emph{prime tuple}. Letting $n$ range over the natural numbers, we wish to see how often \eqref{1.5} is a prime tuple. For instance, consider  $\mathcal{H}=\{0,1\}$ and the tuple $(n,n+1)$.  If $n=2$, we have the prime tuple $(2\, ,3)$. Notice that this is the only prime tuple of this form because, for $n>2$,  one of the numbers $n$ or $n+1$ is an even number bigger than 2.  On the other hand, if $\mathcal{H}= \{0,2\}$,  then we expect that there are infinitely many prime tuples of the form $(n,n+2)$.  This is the twin prime conjecture.
In general, the tuple \eqref{1.5} can be a prime tuple for more than one $n$ only if for every prime $p$ the $h_i$'s  never occupy all of the residue classes modulo $p$.  This is immediately true for all primes $p>k$, so to test this condition we need only to examine small primes. If we denote by 
$ \nu_p(\mathcal{H})$  the number of distinct residue classes modulo $p$ occupied by the integers $h_i$, then we can avoid $p$  dividing some  component of \eqref{1.5} for every $n$ by requiring
\begin{equation}
\nu_p(\mathcal{H}) <p \ \ \text{for all primes}\ p. \label{1.6}
\end{equation}
If this condition holds we say that $\mathcal{H}$ is \emph{admissible} and we call the tuple \eqref{1.5} corresponding to this $\mathcal{H}$ an \emph{admissible tuple}. It is a long-standing conjecture that admissible tuples will infinitely often be prime tuples. Our first result is a step towards confirming this conjecture. 
\begin{theorem}
Suppose the primes have level of distribution $\vartheta > 1/2$.
Then there exists an explicitly calculable constant $C(\vartheta)$
depending only on $\vartheta$ such that any 
admissible $k$-tuple with $k \geq C(\vartheta)$ contains at least two primes infinitely often. 
Specifically, if  $\vartheta \ge 0.971$, then this is true for $k \ge 6$. 
\end{theorem}
Since the 6-tuple $(n,n+4,n+6, n+10, n+12, n+16)$ is admissible, the Elliott-Halberstam conjecture implies that 
\begin{equation}
\liminf_{n\to \infty} (p_{n+1} - p_n) \leq 16,
\label{1.7}
\end{equation}
where the notation $p_n$ is used to denote the $n$-th prime.
This means that $p_{n+1} - p_n \leq 16$ for infinitely many $n$.
Unconditionally, we prove a long-standing conjecture concerning gaps between consecutive primes.
\begin{theorem}  We have 
\begin{equation}E_1 := \liminf_{n\to \infty} \frac{p_{n+1}-p_n}{\log p_n} =0. \label{1.8} \end{equation}
\end{theorem}
There is a long history of results on this topic which we will briefly mention. The inequality $E_1 \leq 1$ is a trivial consequence of the 
prime number theorem. 
The first result of type $E_1 < 1$ was proved in 
1926 by Hardy and Littlewood \cite{HL-0}, who on assuming the Generalized Riemann 
Hypothesis (GRH) obtained
$E_1 \leq 2/3$. 
This result was improved by Rankin \cite{Rankin} to 
$
E_1 \leq 3/5, $
also assuming the GRH. 
The first unconditional estimate was proved by Erd\H os 
\cite{Erdos} in 1940. Using Brun's sieve, he showed that
$ E_1 < 1 - c $ 
with an unspecified positive explicitly calculable constant~$c$. 
His estimate was improved by Ricci \cite{Ricci} in 1954 to 
$E_1 \leq 15/16.$ 
In 1965 Bombieri and Davenport \cite{BD} refined and made unconditional the method of Hardy and Littlewood by substituting  the Bombieri--Vinogradov theorem  for the GRH, and obtained $E_1\le 1/2$. They also combined their method with the method of Erd\H os and obtained 
$E_1 \leq 0.4665\dots\, $. 
Their result was further refined by Pilt'ai \cite{Pi} to $E_1 \leq 0.4571\dots\, $, Uchiyama \cite{Uc} to $E_1 \leq 0.4542\dots\,$ and in several steps by Huxley \cite{Hu1} \cite{Hu2} to 
yield  
$E_1 \leq 0.4425\dots\, $, 
and finally in 1984 to $E_1 \le .4393\ldots \, $ \cite{Hu3}.
In 1988 Maier \cite{Ma} used his matrix-method to improve 
Huxley's result to  
$E_1 \leq e^{-\gamma} \cdot 0.4425\dots = 0.2484\dots\, $,
where $\gamma$ is Euler's constant. 
Maier's method by itself gives $E_1 \leq e^{-\gamma} = 0.5614\dots\,$. 
The recent version of the method of Goldston and Y\i ld\i r\i m
\cite{GYIII} led, without combination with other methods, to 
$E_1 \leq 1/4 $.

 In a later paper in this series we will prove the quantitative result that 
\begin{equation}
\liminf_{n\to \infty} \frac{p_{n+1}-p_n}{(\log
p_n)^{\frac{1}{2}}(\log\log p_n)^{2}} < \infty.
\label{1.9}
\end{equation}

While Theorem 1 is a striking new result, it also reflects the limitations of our current method. Whether these limitations are real or can be overcome is a critical issue for further investigation. We highlight the following four questions. 

\smallskip\noindent\textbf{Question 1.} Can it be proved
unconditionally by the current method that there are infinitely
often bounded gaps between primes?  Theorem 1 would appear to
be within a hair's breadth of obtaining this result. However, any
improvement in the level of distribution $\vartheta$ beyond $1/2$
probably lies very deep, and even the GRH does not help. Still, there
are  stronger versions of the Bombieri-Vinogradov theorem, as found
in \cite{BFI}, and the circle of ideas used to prove these results,
which may help to obtain this
result.

\smallskip\noindent\textbf{Question 2.} Is $\vartheta= 1/2$ a
true barrier for obtaining primes in tuples? Soundararajan \cite{So}
has demonstrated this is the case for the current argument, but
perhaps more efficient arguments may be devised.

\smallskip\noindent\textbf{Question 3.} Assuming the
Elliott-Halberstam conjecture, can it be proved that there are three
or more  primes in admissible $k$-tuples with large enough $k$? Even
under the strongest assumptions, our method fails to prove anything
about more than two primes in a given tuple.

\smallskip\noindent\textbf{Question 4.} Assuming the
Elliott-Halberstam conjecture, can the twin prime conjecture be
proved using the current approximations?

\smallskip The limitation of our method, identified in Question
3, is the reason we are less successful in finding more than two
primes close together. However, we are able to improve on earlier
results, in particular the recent results in \cite{GYIII}. For
$r\geq 1$, let
\begin{equation}
E_r = \liminf_{n \to \infty} \frac{p_{n+r} - p_n}{\log p_n}.
\label{1.10}
\end{equation}
Bombieri and Davenport \cite{BD} showed $E_r \leq r - 1/2$. This
bound was later improved by Huxley \cite{Hu1,Hu2} to $E_r \leq r -
5/8 + o(1/r)$, by Goldston and Y\i ld\i r\i m \cite{GYIII} to $E_r
\leq (\sqrt{r} - 1/2)^2$, and by  Maier \cite{Ma} to $E_r \leq
e^{-\gamma} \left(r - 5/8 + o \left(1/r\right) \right)$. In proving
Theorem 2 we will also show, assuming the primes have level of
distribution $\vartheta$,
\begin{equation}
E_r \leq \max (r - 2 \vartheta, 0), \label{1.11}
\end{equation} 
and hence unconditionally $E_r \le r-1$.  However, by a more complicated argument, we will prove the following result.
\begin{theorem}
Suppose the primes have level of distribution $\vartheta$. Then for $r \geq 2$,
\begin{equation}
E_r \leq (\sqrt{r} - \sqrt{2\vartheta})^2. \label{1.12}
\end{equation}
In particular, we have unconditionally, for $r \geq 1$,
\begin{equation}
E_r \leq (\sqrt{r} - 1)^2. \label{1.13}
\end{equation} 
\end{theorem}
From \eqref{1.11} or \eqref{1.12} we see that the Elliott-Halberstam conjecture implies that
\begin{equation}E_2 = \liminf_{n \to \infty} \frac{p_{n+2} - p_n}{\log p_n} = 0.
\label{1.14}
\end{equation}
We note that if we couple the ideas of the present work
with Maier's matrix method \cite{Ma} we then expect 
that  \eqref{1.12} can be replaced by the stronger inequality
\begin{equation}
E_r \leq e^{-\gamma} (\sqrt{r} - \sqrt{2\vartheta} )^2.
\label{1.15}
\end{equation}

While this paper is our first paper on this subject, we have two other papers which overlaps it.  The  first paper \cite{GMPY},  written jointly with  Motohashi, gives a short and simplified proof of Theorems 1 and 2.  The second paper \cite{GGPY}, written jointly with Graham, uses sieve methods to prove Theorems 1 and 2 and provides applications for tuples of almost-primes (products of two or more distinct primes.) 

The present paper is organized as follows. In Section 2, we describe our method and its relation to earlier work. We also state Propositions 1 and 2 which incorporate the key new ideas in this paper. These are developed in a more general form than in \cite{GGPY} or \cite{GMPY} so as to be employable in many applications.  In Section 3,  we prove Theorems 1 and 2 using these propositions. The method of proof is due to Granville and Soundararajan. In Section 4 we make some further comments on the method used in Section 3. In Section 5 we prove two lemmas needed later. In Section 6, we prove a special case of Proposition 1 which illustrates the key points in the general case. In Section 7 we begin the proof of Proposition 1 which is reduced to evaluating a certain contour integral.
In Section 8 we evaluate a more general contour integral that occurs in the proof of both propositions. In Section 9, we prove Proposition 2.  In this paper we do not obtain results that are uniform in $k$, and therefore we assume here that our tuples have a fixed length. However, uniform results are needed for \eqref{1.9}, and they will be the topic of the next paper in this series. 
Finally, we prove Theorem 3 in Section 10. 

\emph{Notation.} In the following $c$ and $C$ will denote (sufficiently) small and
(sufficiently) large absolute positive constants, respectively, which have been chosen appropriately. This is also true for constants formed from $c$ or $C$ with subscripts or accents.  We unconventionally will allow these constants to be different at different occurences.
Constants implied by pure $o$, $O$, $\ll$ symbols will be absolute, unless otherwise stated. $[S]$ is 1 if the statement $S$\marginpar{} is true and is 0 if
$S$ is false.
 The symbol $\sum^\flat$ indicates the summation is over  squarefree
 integers, and $\sum'$ indicates the summation variables are pairwise relatively prime.

The ideas used in this paper have developed over many years.  We are indebted to many people, not all of whom we can mention. However, we would like to thank A. Balog, E. Bombieri, T. H. Chan, J. B. Conrey, P. Deift, D. Farmer, K. Ford, J. Friedlander, A. Granville, C. Hughes, D. R. Heath-Brown, A. Ledoan, H. L. Montgomery, Sz. Gy. Revesz, P. Sarnak, and K. Soundararajan.

\section{Approximating prime tuples}

Let
\begin{equation}
\mathcal{H}=\{h_1,h_2, \ldots, h_k\} \ \ \textrm{with}\ 1\le
h_1,h_2, \dots, h_k\le h \ \text{distinct integers}, \label{2.1}
\end{equation}
and let $\nu_p(\mathcal{H})$ denote the number of distinct
residue classes modulo $p$ occupied by the elements of
$\mathcal{H}$.\footnote{ The restriction of the set $\mathcal{H}$ to positive integers is only for simplicity, and, if desired, can easily be removed later from all of our results.} For squarefree integers $d$, we extend this
definition to $\nu_d(\mathcal{H})$ by multiplicativity. We denote by
\begin{equation}
\gs(\mathcal{H}) := \prod_p\biggl(1-\frac{1}{p}\biggr)^{-k}\biggl(1
- \frac{\nu_p(\mathcal{H})}{p}\biggr) \label{2.2}
\end{equation}
 the singular series associated with $\mathcal{H}$.
Since $\nu_p(\mathcal{H})=k$ for $p>h$, we see that the product  is absolutely convergent and therefore   $\mathcal{H}$ is admissible as defined in \eqref{1.6} if and only if $\gs(\mathcal{H})\neq 0$.
Hardy and Littlewood conjectured an asymptotic formula for the number of prime tuples $(n+h_1,n+h_2, \ldots , n+h_k)$, with
$1\le n\le N$, as $N\to \infty$.  Let $\Lambda(n)$ denote the von Mangoldt
function which equals $\log p$ if $n=p^m$, $m\geq 1$, and zero otherwise. We define
\begin{equation}
\Lambda(n;\mathcal{H}) := \Lambda(n+h_1)\Lambda(n+h_2)\cdots
\Lambda(n+h_k) \label{2.3}
\end{equation}
and use this function to detect prime tuples and tuples with prime
powers in components, the latter of which can be  removed in applications.
The Hardy--Littlewood prime-tuple conjecture \cite{HL} can be stated in the form
\begin{equation}
\sum_{n\le N}\Lambda(n;\mathcal{H}) = N (\gs(\mathcal{H})+o(1)),
\quad  \mbox{as \ $N\to \infty$.}\label{2.4}
\end{equation}
(This conjecture is trivially true if $\mathcal{H}$ is not admissible.)
Except for the prime number theorem (1-tuples), this conjecture
is unproved.\footnote{ Asymptotic results for the number of primes in tuples, unlike
the existence result in Theorem 1, are beyond the reach of our method.}

The program the first and third authors
have been working on since 1999 is to compute
approximations for \eqref{2.3} using short divisor sums and
apply the results to problems on primes. The simplest
approximation of $\Lambda(n)$ is based on the elementary formula
\begin{equation} \Lambda(n) = \sum_{d|n} \mu(d) \log \frac{n}{d},
\label{2.5}\end{equation}
which can be approximated with the  smoothly truncated divisor sum
\begin{equation} \Lambda_R(n) = \sum_{\substack{d|n\\ d\le R}
}\mu(d) \log \frac{R}{d} .\label{2.6} \end{equation}
Thus, an approximation for $\Lambda(n;\mathcal{H})$ is given by
\begin{equation}  \Lambda_R(n+h_1)\Lambda_R(n+h_2)\cdots
\Lambda_R(n+h_k). \label{2.7}
\end{equation}
 In \cite{GYIII}, Goldston and Y\i ld\i r\i m applied \eqref{2.7}
to detect small gaps between primes and proved
\[ E_1=\liminf_{n\to \infty} \left(\frac{p_{n+1}-p_n}{\log p_n}\right) \le \frac{1}{4}.\]
In the process of that work, they realized that for some applications there might be much
better approximations for prime tuples than \eqref{2.7}, but the
approximation they devised was unsuccessful.  Recently,
 the current authors were able to obtain such an approximation, which is applied here to the problem of small gaps between
primes.

The idea for our new approximation came from a paper of
Heath-Brown \cite{HB} on almost prime tuples.
His result is  itself  a generalization of Selberg's
proof from 1951 (see \cite{Se}, p. 233--245) that the polynomial $n(n+2)$ will infinitely often have at
most five distinct prime factors, so that the same is true for the tuple
 $(n,n+2)$. Not only does our approximation have its origin in these papers, 
but in hindsight the argument of Granville and Soundararajan (employed in the proof of Theorems 1 and 2) is essentially the same as the method used in these papers. 

In connection with the tuple 
\eqref{1.5}, we consider the polynomial
\begin{equation} P_{\mathcal{H}}(n) = (n+h_1)(n+h_2) \cdots
(n+h_k).\label{2.8}
\end{equation}
If the tuple \eqref{1.5} is a prime tuple then $P_{\mathcal{H}}(n)$ has exactly $k$ prime factors.
We detect this condition by using the $k$-th
generalized von Mangoldt function
\begin{equation}
\Lambda_{k}(n) = \sum_{d|n} \mu(d) \biggl(\log
\frac{n}{d} \biggr)^{k },
\label{2.9}
\end{equation}
which vanishes if $n$ has more than $k$ distinct prime factors.\footnote{As with $\Lambda(n)$, we overcount the prime tuples by including factors which are proper prime powers, but these can be removed in applications with a negligible error. The slightly misleading notational conflict between the generalized von Mangoldt function $\Lambda_k$ and $\Lambda_R$ will only occur in this section.}
With this,  our prime tuple detecting function becomes
\begin{equation}
\Lambda_k(n; \mathcal{H}) := \frac{1}{k!} \Lambda_{k } (P_{\mathcal{H}}(n)).
\label{2.10}
\end{equation}
The normalization factor $1/k!$ simplifies the
statement of our results.
As we will see in Section 5,
this approximation suggests the  Hardy--Littlewood type conjecture
\begin{equation}
\sum_{n\le N} \Lambda_k(n;\mathcal{H}) =  N
\left(\mathfrak{S}(\mathcal{H}) +o(1) \right).
\label{2.11}
\end{equation}
This is a special case of the general conjecture of Bateman--Horn \cite{BH} which is the quantitative form of Schinzel's conjecture \cite{Sc}.
There is not much difference between \eqref{2.4}  and \eqref{2.11}, but the same is not true of their approximations. 

In analogy with  \eqref{2.6} (when $k=1$),   we
approximate $\Lambda_k$  by the smoothed and truncated divisor sum
\[ \sum_{\substack{d|n\\ d\le R} }\mu(d) \left(\log \frac{R}{d}
\right)^{k }
\]
and  define
\begin{equation}
\Lambda_R(n; \mathcal{H})
= \frac{1}{k !}\sum_{\substack{ d|P_{\mathcal{H}}(n)\\ d\le R}}
 \mu(d) \left(\log \frac{R}{d}\right)^{k } .
\label{2.12}
\end{equation}
However, as we will see in the next section, this approximation is not adequate to prove Theorems 1 and 2. 

A second simple but crucial idea is needed: rather than only approximate prime tuples, one should  approximate tuples with primes in many components. Thus, we consider when  $P_{\mathcal{H}}(n)$ has $k+\ell$ or less distinct prime factors,  where $0\le \ell \le  k $, and define 
\begin{equation}
\Lambda_R(n; \mathcal{H}, \ell)
= \frac{1}{(k + \ell)!}\sum_{\substack{ d|P_{\mathcal{H}}(n)\\ d\le R}}
 \mu(d) \left(\log \frac{R}{d}\right)^{k + \ell} ,
\label{2.13}
\end{equation}
where $|\mathcal H|=k$. 
In Section 4, we will give precisely a measure of how well a function  detects primes in tuples, so that in terms of this measure, when $k,\ell \to \infty$ and $\ell = o(k)$, this approximation is twice as good at detecting primes in tuples as \eqref{2.12} (that is, when $\ell=0$), which in turn is twice as good as \eqref{2.7}. This improvement enables us to prove Theorem 2 unconditionally. Moreover, it allows the level of distribution needed in Theorem 1 to be any number $>1/2$. 
 
The advantage of \eqref{2.12} and \eqref{2.13} over \eqref{2.7} can be seen as follows.  If  in \eqref{2.12} or \eqref{2.13} we restrict ourselves
to $d$'s with all prime factors larger than $h$, then the condition
$d|P_{\mathcal{H}}(n)$ implies that
we can write $d = d_1 d_2 \cdots d_k$ uniquely with $d_i|
n+h_i$, $1\le i\le k$, the $d_i$'s pairwise relatively prime,
and $d_1d_2 \cdots d_k \le R$.  In our application to prime gaps
we require that $R\le N^{\frac{1}{4} - \epsilon}$. On the other hand,
 on expanding, \eqref{2.7} becomes
a sum over $d_i| n+h_i$, $1\le i\le k$, with
$d_1\le R$, $d_2\le R$, $\ldots $ , $d_k\le R$.  The application
to prime gaps  here requires that $R^k \le
N^{\frac{1}{4} - \epsilon}$, and so $R\le
N^{\frac{1}{4k} - \frac{\epsilon}{k}}$. Thus \eqref{2.7} has a more severe restriction on the range of the
divisors. An additional technical advantage is that having one
truncation rather than $k$ truncations simplifies our calculations.

Our main results on $\Lambda_R(n; \mathcal{H}, \ell)$ are summarized in
the following two propositions. Suppose $\mathcal{H}_1$ and
$\mathcal{H}_2$ are, respectively,  sets of
$k_1$ and $k_2$ distinct non-negative integers $\le h$.  We always assume that at least one of these sets
is non-empty.
Let $M = k_1 + k_2 + \ell_1 + \ell_2$.
\begin{proposition} 
Let $\mathcal{H}=\mathcal{H}_1 \cup \mathcal{H}_2$, $|\mathcal H_i| 
= k_i$, and $r=|\mathcal{H}_1\cap \mathcal{H}_2|$. If $R\ll 
N^{\frac{1}{2}}(\log N)^{-4M}$ and $h\le R^C$ for any given constant 
$C > 0$, then as $R,N\to \infty$ we have
\begin{equation} 
\sum_{n\le N} 
\Lambda_R( n;\mathcal{H}_1,\ell_1)\Lambda_R( n;\mathcal{H}_2,\ell_2) = 
{\ell_1 + \ell_2 \choose \ell_1} \frac{(\log R)^{r + \ell_1 + 
\ell_2}}{(r + \ell_1 + \ell_2)!} (\mathfrak S(\mathcal H) + o_M(1))N 
.\label{2.14} 
\end{equation} 
\end{proposition} 
\begin{proposition} 
Let $\mathcal{H}=\mathcal{H}_1 \cup \mathcal{H}_2$, $|\mathcal{H}_i| 
= k_i$,  $r=|\mathcal{H}_1\cap \mathcal{H}_2|$, $1\le h_0\le h$, and 
$\mathcal{H}^0 = \mathcal{H}\, \cup \, \{h_0\}$.  If $R\ll_M 
N^{\frac{1}{4}}(\log N)^{-B(M)}$ for a sufficiently large positive 
constant $B(M)$, and $h\le R$, then as $R,N\to \infty$ we have 
\begin{equation} \begin{split} 
& \sum_{n\le N} 
\Lambda_R 
( n;\mathcal{H}_1,\ell_1)\Lambda_R( n;\mathcal{H}_2,\ell_2) 
\theta (n+h_0) \\ 
&\quad= \left\{ 
\begin{array}{@{\hspace*{0pt}}l@{\hspace*{2pt}}l@{\hspace*{0pt}}} 
{\displaystyle {\ell_1 + \ell_2 \choose \ell_1} \frac{(\log R)^{r + 
\ell_1 + \ell_2}}{(r + \ell_1 + \ell_2)!} (\mathfrak S(\mathcal H^0) 
+ o_M(1))N} &\quad\mbox{if $h_0\not \in \mathcal{H}$,}\\ 
 {\displaystyle 
{\ell_1 \! +\!  \ell_2 \! +\!  1 \choose \ell_1 + 1} \frac{(\log 
R)^{r + \ell_1 + \ell_2 + 1}}{(r \! +\!  \ell_1 \! +\!  \ell_2 \! 
+\!  1)!} (\mathfrak S(\mathcal H) \! +\!  o_M(1))N} 
 &\quad\mbox{if $h_0 \in \mathcal{H}_1$ and $h_0\not \in 
 \mathcal{H}_2$,}\\ 
 {\displaystyle 
{\ell_1 \! +\!  \ell_2 \! +\!  2 \choose \ell_1 + 1} \frac{(\log 
R)^{r + \ell_1 + \ell_2 + 1}}{(r \! +\!  \ell_1 \! +\!  \ell_2 \! 
+\!  1)!} (\mathfrak S(\mathcal H) \! +\!  o_M(1))N} 
 &\quad\mbox{if $h_0 \in \mathcal{H}_1\cap \mathcal{H}_2$.}\\ 
\end{array} 
\right.\end{split}\label{2.15} 
\end{equation} 
Assuming the primes have level of distribution $\vartheta >1/2$, 
i.e., \eqref{1.3} with \eqref{1.4} holds, we may choose, for any 
$\epsilon>0$, $R\ll_M N^{\frac{\vartheta}{2} - \epsilon}$ and $h\le R^\epsilon$.
\end{proposition} 

\begin{rema}
By relabeling the variables we obtain the corresponding form
if $h_0 \in \mathcal H_2, h_0 \not\in \mathcal H_1$.
\end{rema}

Propositions 1 and 2 can be strengthened in several ways. We will show that the error terms $o_M(1)$ can be replaced by a series of lower order terms and a prime number theorem type of error term. Moreover, we can make the result uniform for $M\to \infty$ as an explicit function of $N$ and $R$. This will be proved in a later paper and used in the proof of \eqref{1.9}.

\section{Proofs of Theorems 1 and 2}

In this section we employ Propositions 1 and 2 and
a simple argument due to Granville and Soundararajan 
 to prove Theorems 1 and 2.

For $\ell \geq 0$, $\mathcal H_k = \{h_1, h_2, \dots, h_k\}$, $1\le h_1,h_2,\ldots , h_k\le h\le R$,
we deduce from
Proposition~1, for $R \ll N^{\frac{1}{2}} (\log N)^{-B(M)}$ and $R,N\to \infty$, that
\begin{equation}
\sum_{n\leq N} \Lambda_{R} (n; \mathcal H_k, \ell)^2 \sim
\frac1{(k + 2\ell)!} {2\ell\choose \ell} \mathfrak S(\mathcal
H_k) N (\log R)^{k + 2\ell}.
\label{3.1}
\end{equation}
For any  $h_i \in \mathcal H_k$, we have from Proposition~2, for $R\ll N^{\frac{\vartheta}{2} -\epsilon}$, and $R,N\to \infty$, 
\begin{equation}
\sum_{n\leq N} \Lambda_R(n; \mathcal H_k, \ell)^2 \theta (n +
h_i) \sim \frac2{(k + 2\ell + 1)!} {2\ell + 1\choose \ell}
\mathfrak S(\mathcal H_k) N(\log R)^{k+2\ell+1}.
\label{3.2}
\end{equation}
Taking $R= N^{\frac{\vartheta}{2} -\epsilon}$, we obtain\footnote{In \eqref{3.3}, as well as later in \eqref{3.8}, the asymptotic sign replaces an error term of size $o(\log N)$ in the parenthesis term after $\log 3N$. We thus make the convention that the asymptotic relationship holds only up to the size of the apparent main term.}
\begin{equation}\begin{split}
\mathcal S:&= \sum^{2N}_{n = N+1} \left( \sum^k_{i = 1}
\theta
(n + h_i) - \log 3N\right) \Lambda_R(n; \mathcal H_k, \ell)^2 
\\
& \sim k \frac2{(k + 2\ell + 1)!} {2\ell + 1\choose \ell}
\mathfrak S(\mathcal H_k) N(\log R)^{k+2\ell+1}\\
&\quad - \log 3N
\frac{1}{(k + 2\ell)!} {2\ell \choose \ell} \mathfrak S(\mathcal
H_k) N(\log R)^{k+2\ell} \\
& \sim \left(\frac{2k}{k + 2\ell + 1} \frac{2\ell + 1}{\ell + 1}
\log R - \log 3N\right) \frac1{(k + 2\ell)!} {2\ell \choose
\ell} \mathfrak S(\mathcal H_k) N(\log R)^{k + 2\ell}. \label{3.3}
\end{split}\end{equation}
Here we note that the tuple $\mathcal H_k$ will contain at least two primes if
$\mathcal S > 0$. This situation occurs when
\begin{equation}
\frac{k}{k+2\ell + 1} \frac{2\ell + 1}{\ell + 1} \vartheta > 1.
\label{3.4}
\end{equation}
If $k, \ell \to \infty$ with $\ell = o(k)$, then the left-hand side has the limit $2\vartheta$, and thus \eqref{3.4} holds for any $\vartheta > 1/2$ if we choose $k$ and $\ell$ appropriately depending on $\vartheta$.  This proves the first part of Theorem 1.
Next, assuming $\vartheta > 20/21$, we see that \eqref{3.4} holds
with $\ell = 1$ and $k = 7$. This proves the second part of Theorem
1 but with $k=7$. The case $k=6$ requires a slightly more
complicated argument and is treated later in this section.

The table below gives the values of $C(\vartheta)$, defined in Theorem 1, obtained from \eqref{3.4}.
For a given $\vartheta$,  it gives the smallest $k$ and corresponding smallest $\ell$ for which \eqref{3.4} is true. Here $h(k)$ is the shortest length of any admissible $k$-tuple,
which has been computed by Engelsma \cite{En} by exhaustive search
for $1\le k \le 305$ and covers every value in this table and the
next except $h(421)$, where we have taken the upper bound value from
\cite{En}.
\smallskip 
\begin{center}
\begin{tabular}{|l|l|l|l|}\hline
$\vartheta$ & $k$ & $\ell$& $h(k)$\\
\hline \hline
1 &7&1&20\\
\hline
.95&8&1&26\\
\hline
.90&9&1&30\\
\hline
.85&11&1&36\\
\hline
.80&16&1&60\\
\hline
.75&21&2&84\\
\hline
.70&31&2&140\\
\hline
.65&51&3&252\\
\hline
.60&111&5&634\\
\hline
.55&421&10&$2956^*$\\
\hline
\end{tabular}
\smallskip  
\end{center}
\begin{center}* indicates this value could be an upper bound of the true value. \end{center}
\medskip
To prove Theorem 2, we modify the previous proof by considering
\begin{equation}
\widetilde{\mathcal{S}}
:= \sum_{n=N+1}^{2N}
\left( \sum_{1\le h_0\le h} \theta (n+h_0) - r \log 3N \right)
\sum_{\substack{1\le h_1,h_2, \ldots ,h_k\le h \\
\text{distinct}}}
\Lambda_{R}(n;\mathcal{H}_k, \ell)^2 ,
\label{3.5}
\end{equation}
where $r$ is a positive integer. 
 To evaluate $\widetilde{\mathcal{S}}$, we need the case of Proposition 2 where
$h_0\not \in \mathcal{H}_k$, 
\begin{equation}
\sum_{n\leq N} \Lambda_R(n; \mathcal H_k, \ell)^2\,  \theta (n +
h_0) \sim \frac1{(k + 2\ell )!} {2\ell \choose \ell}
\mathfrak S(\mathcal{H}_k\cup \{h_0\}) N(\log R)^{k+2\ell}.
\label{3.6}
\end{equation}
We also need a result of Gallagher \cite{Ga}: as $h\to \infty$,
\begin{equation}
\sum_{\substack{1\le h_1,h_2, \ldots h_k\le h \\
\text{distinct}}}
\mathfrak{S}(\mathcal{H}_k ) \sim h^k.
\label{3.7}
\end{equation}
Taking $R= N^{\frac{\vartheta}{2} - \epsilon}$, and applying \eqref{3.1}, \eqref{3.2}, \eqref{3.6}, and \eqref{3.7}, we find that
\begin{equation}\begin{split} 
 \widetilde{\mathcal{S}}&\sim \sum_{\substack{1\le h_1,h_2, \ldots 
,h_k\le h \\ \text{distinct}}} 
\Bigg(k\frac{2}{(k+2\ell+1)!}\genfrac{(}{)}{0pt 
}{0}{2\ell+1}{\ell} \mathfrak{S}(\mathcal{H}_k) N(\log 
R)^{k+2\ell+1} 
\\ 
& \quad +\sum_{\substack{1\le h_0\le h \\ h_0\neq h_i, 1\le i\le k 
}} \frac{1}{(k+2\ell)!}\genfrac{(}{)}{0pt }{0}{2\ell}{\ell} 
\mathfrak{S}(\mathcal{H}_k\cup \{h_0\}) N(\log R)^{k+2\ell} \\ 
& \qquad  - r\log 3N  \frac{1}{(k+2\ell)!}\genfrac{(}{)}{0pt 
}{0}{2\ell}{\ell} \mathfrak{S}(\mathcal{H}_k) N(\log 
R)^{k+2\ell}\Bigg)\\ 
& \sim  \left( \frac{2k}{k+2\ell+1}\frac{2\ell+1}{\ell +1} 
\log R + h - r\log 
3N\right)\frac{1}{(k+2\ell)!}\genfrac{(}{)}{0pt 
}{0}{2\ell}{\ell} Nh^k(\log R)^{k+2\ell}. \label{3.8} 
\end{split}\end{equation}
 Thus, there
are at least $r+1$ primes in some interval $(n,n+h]$, $N<n\le 2N$, provided that
\begin{equation}
h  > \biggl(r -\frac{2k}{k+2\ell+1}\frac{2\ell+1}{\ell +1}
\biggl(\frac{\vartheta}{2}-\epsilon\biggr) \biggr)\log N,
\label{3.9}
\end{equation}
which, on letting $\ell = [\sqrt{k}/2]$ and taking $k$ sufficiently large,
gives
\begin{equation} h> \left(r -  2\vartheta  +4 \epsilon +O(\frac{1}{\sqrt{k}}) \right)\log N.
\label{3.10}
\end{equation}
This proves \eqref{1.11}.  Theorem 2 is the special case $r=1$ and $\vartheta=1/2$.
 
We are now ready to prove the last part of Theorem 1. Consider
\begin{equation}\begin{split}
\mathcal S':&= \sum^{2N}_{n = N+1} \left( \sum^k_{i = 1}
\theta
(n + h_i) - \log 3N\right) \left( \sum_{\ell =0}^ L a_\ell
\Lambda_R(n; \mathcal H_k, \ell)\right)^2  
\\
&=\sum^{2N}_{n = N+1} \left( \sum^k_{i = 1}
\theta
(n + h_i) - \log 3N\right)  \sum_{0\le \ell_1,\ell_2 \le L}
a_{\ell_1}a_{\ell_2} \Lambda_R(n; \mathcal H_k,
\ell_1)\Lambda_R(n; \mathcal H_k, \ell_2)\\
& = \sum_{0\le \ell_1,\ell_2 \le L} a_{\ell_1}a_{\ell_2}
\mathcal{M}_{\ell_1,\ell_2}, \label{3.11}
\end{split}\end{equation}
where
\begin{equation}
\mathcal{M}_{\ell_1,\ell_2} = \widetilde{M}_{\ell_1,\ell_2} -
(\log 3N) M_{\ell_1,\ell_2} ,
\label{3.12} \end{equation}
say.
Applying Propositions 1 and 2 with $R=
N^{\frac{\vartheta}{2}-\epsilon}$, we
deduce that
\[ M_{\ell_1,\ell_2} \sim
{\ell_1 + \ell_2 \choose \ell_1} \frac{(\log R)^{k + \ell_1 +
\ell_2}}{(k + \ell_1 + \ell_2)!} \mathfrak S(\mathcal{ H}_k)N\]
and
\[ \widetilde{M}_{\ell_1,\ell_2} \sim
k{\ell_1 + \ell_2 +2\choose \ell_1+1} \frac{(\log R)^{k + \ell_1 +
\ell_2+1}}{(k + \ell_1 + \ell_2+1)!} \mathfrak S(\mathcal{ H}_k)N.\]
Therefore,
\[ \begin{split}\mathcal{M}_{\ell_1,\ell_2}
 \sim {\ell_1 +
\ell_2 \choose \ell_1} &\mathfrak{S}(\mathcal{H}_k)N\frac{(\log
R)^{k + \ell_1 +
\ell_2}}{(k + \ell_1 + \ell_2)!}\\& \times\left(\frac{k( \ell_1
+ \ell_2+2)(\ell_1 + \ell_2+1)}{(\ell_1+1)(\ell_2+1)(k + \ell_1
+ \ell_2+1)}\log R -\log 3N\right).\end{split}\]
Defining $b_\ell = (\log R)^\ell a_\ell$ and $\boldsymbol{b}$ to be the column matrix corresponding to the vector 
$(b_0,b_1,\ldots, b_L)$,  we obtain
\begin{equation}
\begin{split}
&S^*(N,\mathcal{H}_k,\vartheta, \boldsymbol{b})
:=\frac{1}{\mathfrak{S}(\mathcal{H}_k)N(\log R)^{k+1}} S'\\
&\hspace{0.1cm} \sim \sum_{0\le \ell_1,\ell_2\le L}\!\!
b_{\ell_1}b_{\ell_2}{\ell_1+\ell_2 \choose \ell_1}
\frac{1}{(k+\ell_1+\ell_2)!}\left(\frac{k(\ell_1+\ell_2+2)(\ell_1 +
\ell_2+1)}{(\ell_1 + 1)(\ell_2+1)(k+\ell_1+\ell_2+1)} -
\frac{2}{\vartheta}\right)\\
&\hspace{0.1cm} \sim \boldsymbol{b}^T \boldsymbol{M} \boldsymbol{b},
\label{3.13}
\end{split}
\end{equation}
where
\begin{equation} \boldsymbol{M} = \left[ {i+j \choose i}
\frac{1}{(k+i+j)!}\left(\frac{k(i+j+2)(i+j+1)}{(i+1)(j+1)(k+i+j+1)}
- \frac{2}{\vartheta}\right)\right]_{
0\le i,j\le L}.\label{3.14}\end{equation}
We need to choose $\boldsymbol{b}$ so that $S^*>0$ for a given
$\vartheta$ and minimal $k$. On taking $\boldsymbol{b}$ to be an eigenvector of the matrix
$\boldsymbol{M}$ with eigenvalue $\lambda$, we see that
\begin{equation}
S^* \sim \boldsymbol{b}^T \lambda \boldsymbol{b}  = \lambda
\sum_{i=0}^k |b_i|^2
\label{3.15} \end{equation}
 will be $>0$ provided that $\lambda$ is positive. Therefore
$S^* >0$ if $\boldsymbol{M}$ has a positive eigenvalue and
$\boldsymbol{b}$ is chosen to be the corresponding eigenvector. 
Using \emph{Mathematica} we computed the values of $C(\vartheta)$ indicated in the following table, which may be compared to the 
earlier table obtained from \eqref{3.4}. 
\smallskip
\begin{center}
\begin{tabular}{|l|l|l|l|}\hline
$\vartheta$ & $k$ & $L$& $h(k)$\\
\hline \hline
1 &6&1&16\\
\hline
.95&7&1&20\\
\hline
.90&8&2&26\\
\hline
.85&10&2&32\\
\hline
.80&12&2&42\\
\hline
.75&16&2&60\\
\hline
.70&22&4&90\\
\hline
.65&35&4&158\\
\hline
.60&65&6&336\\
\hline
.55&193&9&1204\\
\hline
\end{tabular}
\end{center}
\vspace{0.3cm} In particular, taking $k=6$, $L=1$, $b_0=1$, and
$b_1=b$ in
\eqref{3.13}, we get
\[ \begin{split} S^* &\sim \frac{1}{8!}\left(96-\frac{112}
{\vartheta} +2b\left(18-\frac{16}{\vartheta}\right)
+b^2\left(4-\frac{4} {\vartheta}\right)\right)\\& \sim
-\frac{4(1-\vartheta)}{8!\vartheta}\left(b^2 -
2b\frac{18\vartheta-16}{4(1-\vartheta)} -\frac{96\vartheta
-112}{4(1-\vartheta)} \right)\\& \sim
-\frac{4(1-\vartheta)}{8!\vartheta}\biggl(\left(b-\frac{18\vartheta-16}{4(1
-\vartheta)}\right)^2 +\frac{15\vartheta^2 -64\vartheta
+48}{4(1-\vartheta)^2}\biggr)
.\end{split}\]
Choosing $b=\frac{18\vartheta-16}{4(1-\vartheta)}$, we then have
\[ S^* \sim - \frac{15\vartheta^2 -64\vartheta
+48}{8!\vartheta(1-\vartheta)},
\]
of which the right-hand side is $>0$ if $\vartheta\le 1$ lies
between the two roots of the quadratic; this occurs when
$4(8-\sqrt{19})/15<\vartheta \le 1$. Thus, there are at least two
primes in any admissible tuple $\mathcal{H}_k$ for $k = 6$, if
\begin{equation}
\vartheta > \frac{4(8-\sqrt{19})}{15}=0.97096\ldots. \label{3.16}
\end{equation}
This completes the proof of Theorem 1. 

\section{Further Remarks on Section 3}

We can formulate the method of Section 3 as follows. For a given tuple $\mathcal{H}=\{h_1,h_2,\ldots , h_k\}$ we define 
\begin{equation}  Q_1 := \sum_{n=N+1}^{2N} f_R(n; \mathcal{H})^2, \quad
Q_2 := \sum_{n=N+1}^{2N}\left(\sum_{i=1}^k\theta(n+h_i)\right) f_R(n; \mathcal{H})^2 ,
\label{4.1}\end{equation}
where $f$ should be chosen to make $Q_2$ large compared with $Q_1$, and $R=R(N)$ will be chosen later. It is reasonable to assume
\begin{equation} f_R(n;\mathcal{H}) = \sum_{\substack{d|P_{\mathcal{H}}(n) \\ d\le R}} \lambda_{d,R} .\label{4.2} \end{equation}
Our goal is to select the $\lambda_{d,R}$ which  maximizes 
\begin{equation} \rho = \rho(N;\mathcal{H},f) := \frac{1}{\log 3N}\left(\frac{Q_2}{ Q_1}\right).\label{4.3} \end{equation}
If $\rho >r$ for some $N$ and  positive integer $r$,  then there exists an $n$, $N<n\le 2N$, such that the tuple \eqref{1.5} has at least $r+1$ prime components.  

This method is exactly the same as that introduced for twin primes by Selberg and for general tuples by Heath-Brown. However,  they used the divisor function $d(n+h_i)$ in $Q_2$ in place of $\theta(n+h_i)$.  Heath-Brown even chose  $f=\Lambda_R(n;\mathcal{H},1)$. 

As a first example, suppose we choose $f$ as in \eqref{2.6} and \eqref{2.7}, so that
\begin{equation} f_R(n; \mathcal{H}) = \prod_{i=1}^k \Lambda_R(n+h_i) .\label{4.4} \end{equation}
By \cite{GYIII}, we have, as $R, N\to \infty$, \footnote{For special reasons, the validity of the formula for $Q_2$  actually holds here for $R\le N^{\frac{\vartheta}{2(k-1)}(1-\epsilon)}$ if $k\ge 2$, but this is insignificant for the present discussion.}
\begin{equation}
\gathered Q_1 \sim N\mathfrak{S}(\mathcal{H})(\log R)^k
\hspace{0.9cm} \text{if} \ R\le N^{\frac{1 }{2k}(1-\epsilon)}, \\
Q_2 \sim kN\mathfrak{S}(\mathcal{H})(\log R)^{k+1} \quad \text{if} \
R\le N^{\frac{\vartheta }{2k}(1-\epsilon)}.
\endgathered
\label{4.5}
\end{equation}
On taking $R= N^{\frac{\vartheta_0 }{2k}}$, $0<\vartheta_0 < \vartheta$, we see that, as $N\to \infty$,
\begin{equation} \rho \sim k\frac{\log R}{\log N} \sim \frac{\vartheta_0}{2}. \label{4.6}\end{equation}
Notice that $\rho <1$, so we fail to detect primes in tuples. In
Section 3, we proved that on choosing $f =
\Lambda_R(n;\mathcal{H},\ell) $, by \eqref{3.1} and \eqref{3.2}, as
$N\to \infty$,
\begin{equation}  \rho \sim \frac{k}{k+2\ell+1} \frac{2\ell+1}{\ell +1}\vartheta_0. \label{4.7} \end{equation}
If $\ell=0$ this gives $\rho \sim \frac{k}{k+1}\vartheta_0$, which, for large $k$, is twice as large as \eqref{4.6}, while \eqref{4.7} gains another factor of two when $\ell \to \infty$ slowly as $k\to \infty$.  This finally shows $\rho >1$ if $\vartheta>1/2$, but just fails if $\vartheta=1/2$.

In \eqref{3.11} we chose
\begin{equation}f_R(n;\mathcal{H})  = \sum_{\ell =0}^ L \frac{b_\ell}{(\log R)^\ell}
\Lambda_R(n; \mathcal H_k, \ell)= \sum_{\substack{d|P_{\mathcal{H}}(n) \\ d\le R}} \mu(d)P\left(\frac{\log (R/d)}{\log R}\right) \label{4.8}\end{equation}
where $P$ is a polynomial with a $k$-th order zero at 0. 
The matrix procedure does not provide a method for analyzing $\rho$ unless $L$ is taken fixed, but the general problem has been solved by Soundararajan \cite{So}. In particular, he showed that $\rho <1$ if $\vartheta=1/2$, so that one can not prove there are bounded gaps between primes using \eqref{4.8}. The exact solution from Soundrarajan's analysis was obtained by a calculus of variations argument by Conrey, which gives, as $N\to \infty$,
\begin{equation} \rho = \frac{k(k-1)}{2\beta}\vartheta_0, \label{4.9}\end{equation}
where $\beta$ is determined as the solution of the equation 
\begin{equation}
\frac{1}{\beta} = \frac{ \int_0^1y^{k-2}q(y)^2
\,dy}{\int_0^1y^{k-1}q'(y)^2 \,dy}\quad \text{with}\ q(y) =
J_{k-2}(2\sqrt{\beta}) - y^{1-\frac{k}{2}}J_{k-2}(2 \sqrt{\beta y})
\label{4.10}
\end{equation}
 where $J_k$ is the Bessel function of the first type.
Using \textit{Mathematica}, one can check that this gives exactly the values of $k$ in the previous table,  which is in agreement with our earlier calculations; but it provides somewhat smaller values of $\vartheta$ for which a given $k$-tuple will contain two primes. Thus, for example, we can replace \eqref{3.16} by the result that every admissible 6-tuple will contain at least two primes if
\begin{equation} \vartheta > .95971 \ldots \ . \label{4.11} \end{equation}

\section{Two Lemmas}
In this section we will prove two lemmas needed for the proof of  Propositions 1 and 2. The conditions on these lemmas have been constructed in order for them to hold uniformly in the given variables. 

The Riemann zeta-function has the Euler product representation,
with $s=\sigma +it$,
\begin{equation} \zeta(s) =
\prod_p\left(1-\frac{1}{p^s}\right)^{-1}, \quad \sigma >1.
\label{5.1}\end{equation}
The zeta-function is analytic except for a simple pole at $s=1$,
where as $s\to 1$
\begin{equation} \zeta(s) = \frac{1}{s-1} + \gamma + O(|s-1|)
.\label{5.2} \end{equation}
(Here $\gamma$ is Euler's constant.)
We need
 standard information concerning the classical zero-free 
region of the Riemann zeta-function. By  Theorem 3.11 and
(3.11.8) in \cite{T},  there exists a small constant
$\overline c>0$, for which we assume $\overline c \leq 10^{-2}$,
such that $\zeta(\sigma +it)\neq 0$ in the region
\begin{equation}
\sigma \ge 1 -\frac{4\overline c}{ \log(|t| + 3)}
\label{5.3}\end{equation}
 for all $t$. Furthermore, we have
\begin{equation}
\gathered
 \zeta(\sigma +it) -
\frac{1}{\sigma -1 +it} \ll \log(|t|+3), \quad
\frac{1}{ \zeta(\sigma +it)} \ll \log(|t|+3),
\\
\frac{\zeta'}{\zeta}(\sigma + it) + \frac1{\sigma - 1+it}
\ll  \log(|t|+3) ,
\endgathered
\label{5.4}
\end{equation}
in this region.
We will fix this $\overline c$ for the rest of the paper
(we could take, for instance, $\overline c = 10^{-2}$, see \cite{Fo}).
 Let $\mathcal{ L}$ denote the contour given by
\begin{equation}
s= -\frac{\overline c}{ \log(|t|+3)} +it.
\label{5.5}\end{equation}

\begin{lemma}
We have, for $R\ge C$, $k\ge 2$,  $B \leq Ck$,
\begin{equation}
\int_{\mathcal{L}}(\log(|s|+3))^B\left|  \frac
{R^s }{s^k} \, ds \right| \ll
C^k_{1} R^{-c_2} +  e^{-\sqrt{\overline c \log R}/2},
\label{5.6}
\end{equation}
where $C_1, c_2$ and the implied constant in $\ll $ depends
only on the constant $C$ in the formulation of the lemma. In addition, if $k \leq c_3
\log R$ with a sufficiently small $c_3$ depending only on $C$, then 
\begin{equation}
\int_{\mathcal L}
(\log(|s| + 3))^B\left|  \frac
{R^s }{s^k} \, ds \right| \ll e^{-\sqrt{\overline c \log R}/2}.
\label{5.7}
\end{equation}
\end{lemma}

\noindent
\emph{Proof.}  The left-hand side of \eqref{5.6} is, with $C_4$ depending on $C$,
\begin{equation}
\aligned
&\ll \int_{0}^\infty R^{\sigma(t)}
\frac{(\log(|t|+4))^B}{(|t|+ \overline c/2)^k}\, dt\\
&\ll \int^{C_4}_0  C^k_1 R^{-c_2} dt + \int^{\omega - 3}_{C_4}
\frac{R^{- \frac{\overline c}{\log(|t| + 3)}}}{t^{3/2}} dt +
\int^\infty_{\omega - 3} t^{-3/2} dt \\
& \ll C^k_1 R^{-c_2} + e^{-\frac{\overline c \log R}{\log
\omega}} + \omega^{-\frac12},
\endaligned
\label{5.8}
\end{equation}
where now $C_1$ is a constant depending on $C$. 
On choosing $\log \omega = \sqrt{\overline {c} \log R}$, the first part of the lemma follows.   The second part is an immediate consequence of the first part.

The next lemma provides some explicit estimates for sums of the
generalized divisor function.
Let $\omega(q)$ denote the number of prime factors of a
squarefree integer $q$.  For any real number $m$, we define
\begin{equation}
d_m(q) = m^{\omega(q)}.
\label{5.9}
\end{equation}
This agrees with the usual definition of the divisor functions when $m$ is a positive integer.  Clearly, $d_m(q) $ is  a monotonically
increasing  function of $m$ (for a fixed $q$), and  for real $m_1$, $m_2$, and $y$, we see  that
\begin{equation}
d_{m_1}(q) d_{m_2}(q) = d_{m_1 m_2} (q), \quad (d_m(q))^y = d_{m^y} (q).
\label{5.10}
\end{equation}
Recall $\sum^\flat$ indicates a sum  over  squarefree
integers. We use the ceiling function  $\lceil y \rceil
:= \min \{n \in \mathbb Z; y \leq n\}$. 
\begin{lemma} We have, for any positive real $m$ and $x\ge 1$
\begin{equation} 
D'(x,m) := \sideset{}{^\flat}\sum_{q\leq x} \frac{d_m(q)}{q} \leq (\lceil m \rceil + \log x)^{\lceil m \rceil} \leq (m
+ 1 + \log x)^{m+1}
\label{5.11}
\end{equation}
and
\begin{equation}
D^*(x,m) := \sideset{}{^\flat}\sum_{q\leq x}d_m(q) \leq x( \lceil
m \rceil + \log x)^{\lceil m \rceil} \leq x(m + 1 + \log x)^{m +
1}.
\label{5.12} \end{equation}
For  $\nu \geq \max( c' \log(K + 1),  1)$,
there is an absolute constant $C'$
depending on $c'$ such that, for $x\ge 1$ and $K\ge 1$, we have
\begin{equation}
\sideset{}{^\flat}\sum_{q\leq x} \frac{(d_{3K}(q))^{1+\frac{1}{\nu}}}{q}
\leq (C' K +
\log x)^{C'K}.
\label{5.13} \end{equation}
\end{lemma}

\noindent
{\it Proof.}
First, we treat the case when $m$ is a positive integer. We prove \eqref{5.11} by induction.
Observe the assertion is true for $m = 1$, that is, when $d_1(q) = 1$ by
definition.
Suppose \eqref{5.11} is proved for $m - 1$. Let us denote the smallest
term in a given product representation of $q$ by $j = j(q) \leq x^{1/m}$.
Then this factor can stand at $m$ places,  and, therefore, with
$q = q' j(q) = q' j$, we have
\[
\aligned
\sideset{}{^\flat}\sum_{q\leq x}\,  \frac{d_m(q)}{q}
&\leq m
\sum^{x^{1/m}}_{j = 1}{}^{^{\scriptstyle\flat}\, } \frac1 j \sideset{}{^\flat}\sum_{q'\leq x/j}
 \frac{d_{m -
1}(q')}{q'} \leq m  (1 + \log x^{\frac{1}{m}})
\left(m - 1 + \log {x} \right)^{m-1}\\
&\leq (m + \log x) (m + \log x)^{m - 1}= (m+\log x)^m.
\endaligned
\]
This completes the induction.  For real $m$, the result holds since $D'(x,m)\le D'(x,\lceil m\rceil )$.
Next, we note that \eqref{5.12} follows from \eqref{5.11} because $D^*(x,m)\le x D'(x,m)$.  
To prove \eqref{5.13}, let  $r := (3K)^{1 + \frac{1}{\nu}} \leq 9 e^{1/c'} K$.  By  \eqref{5.9}, we have
\[
\left(d_{3K}(q)\right)^{1+\frac{1}{\nu}} = d_r(q),
\]
and the result follows by \eqref{5.11} with $C' = 9e^{1/c'} + 1$.

\section{A special case of Proposition~1}
In this section we will prove a special case of Proposition~1  which
illustrates the method without involving the technical complications that appear in the general case. This allows us to set up some notation and obtain estimates for use in the general case. We also obtain the result uniformly in $k$.  

Assume
$\mathcal{H}$ is non-empty (so that $k\ge 1$),  $\ell = 0$, and  $\Lambda_R(n; \mathcal H, 0) = \Lambda_R(n; \mathcal H)$. 
\begin{Prop1simple} Supposing
\begin{equation}
k \ll_{\eta_0}(\log R)^{\frac{1}{2} - \eta_0} \text{ with an arbitrarily small
fixed } \eta_0 > 0,
\label{6.1}
\end{equation}
and $h\le R^C$,
with $C$ any fixed positive number, 
we have
\begin{equation}
\sum_{n=1}^N \Lambda_R(n; \mathcal{H})=
\gs(\mathcal{H}) N + O(Ne^{-c\sqrt{\log R}})
+ O\big(R(2\log R)^{2k}\big) .
\label{6.2}
\end{equation}
\end{Prop1simple}
\noindent
This result motivates the conjecture  \eqref{2.11}.
\medskip

\noindent
\textit{Proof}.  We have
\begin{equation}  \mathcal{S}_R(N;\mathcal{H}) :=\sum_{n=1}^N
\Lambda_R(n;\mathcal{H}) =
\frac{1}{k!}\sum_{d\le R}\mu(d)\left(\log \frac{R}{d}\right)^k
\sum_{\substack{1\le n\le N\\  d|P_{\mathcal{H}}(n)}} 1.
\label{6.3} \end{equation}
If for a prime $p$ we have $p|P_{\mathcal{H}}(n)$, then  among
the solutions $n\equiv -h_i (\text{mod}\, p)$, $1\le i\le k$,
there will be $\nu_p(\mathcal{H})$ distinct solutions modulo
$p$.  For $d$ squarefree we then have by multiplicativity
$\nu_d(\mathcal{H})$ distinct solutions for $n$ modulo $d$ which
satisfy  $d|P_{\mathcal{H}}(n)$, and for each solution one has
$n$ running through a residue class modulo $d$. Hence we see that
\begin{equation} \sum_{\substack{1\le n\le N\\  d|P_{\mathcal{H}}(n)}} 1 =
\nu_d(\mathcal{H}) \left( \frac{N}{d} +O(1)\right) .
\label{6.4}\end{equation}
Trivially  $ \nu_q(\mathcal{H}) \le
k^{\omega(q)} =  d_k(q)$ for squarefree $q$. Therefore,   we conclude that
\begin{equation}
\begin{split} \mathcal{S}_R(N;\mathcal{H})& = N
\left(\frac{1}{k!}\sum_{d\le
R}\frac{\mu(d)\nu_{d}(\mathcal{H})}{d}\left(\log
\frac{R}{d}\right)^k\right)  + O\left(\frac{(\log
R)^k}{k!}\sideset{}{^\flat}\sum_{d\le R}\, \nu_d(\mathcal{H})\right)\\
&= N \mathcal{T}_R(N;\mathcal{H}) + O\big(R (k + \log R)^{2k}\big),
\label{6.5}
\end{split}
\end{equation}
by  Lemma 2.

 Let $(a)$ denote
the contour $s = a + it$, $-\infty <
t<\infty$.
We  apply the formula
\begin{equation}
\frac{1}{ 2\pi i}
\mathop{\int}_{(c)}
\frac{x^s}{ s^{k+1}}\, ds =
\left\{ \begin{array}{ll}
       0
        &\text{if $0<x\le 1$,
    }  \\
\frac{1}{k!}(\log x)^k & \text{if $x\ge 1$,}
\end{array}
\right. 
\label{6.6}
\end{equation}
 for  $c>0$, and have that
\begin{equation}
\mathcal{T}_R(N;\mathcal{H})
= \frac {1}{ 2\pi i} \mathop{\int}_{(1)}
F(s)\frac{R^{s}}{ {s}^{k+1}} ds ,
\label{6.7}
\end{equation}
where, letting $s=\sigma+it$ and assuming $\sigma >0$,
\begin{equation}  F(s) = \sum_{d=1}^\infty \frac{\mu(d)
\nu_{d}(\mathcal{H})}{d^{1+s}}
=\prod_p\Big(1 -
\frac{\nu_p(\mathcal{H})}{p^{1+s}}\Big).
\label{6.8}
\end{equation}
Since $\nu_p(\mathcal{H}) =k$ for all $p>h$,  we
 write
\begin{equation} F(s) =  \frac{G_{\mathcal{H}}(s)}{\zeta(1+s)^k
} ,
\label{6.9}\end{equation}
where  by \eqref{5.1}
\begin{equation}  \begin{split} G_{\mathcal{H}}(s) &= \prod_p\left(1 -
\frac{\nu_p(\mathcal{H})}{p^{1+s}}\right)
\left(1-\frac{1}{p^{1+s}}\right)^{-k}\\&
=\prod_p\left(1 +
\frac{k-\nu_p(\mathcal{H})}{p^{1+s}}+ O_h\Big(\frac{k^2}{p^{2+2\sigma}}\Big)\right),
\label{6.10}\end{split}
\end{equation}
which is analytic and uniformly bounded for $\sigma > - 1/2 + \delta$
for any $\delta > 0$. Also, by \eqref{2.2} we see that
\begin{equation} G_{\mathcal{H}}(0) = \gs(\mathcal{H}).
\label{6.11} \end{equation}

From  \eqref{5.4} and \eqref{6.9},  the function $F(s)$ satisfies the bound 
\begin{equation}
F(s) \ll |G_{\mathcal{H}}(s)|(C \log
(|t|+3))^k .
\label{6.12}
\end{equation}
in the region on and to right of $\mathcal{L} $.
Here $G_{\mathcal{H}}(s)$ is analytic and bounded in this
region, and has a dependence on both $k$ and the size $h$ of
the components of $\mathcal{H}$. 
We note that $\nu_p(\mathcal{H})=k$ not only when $p>h$, but whenever
$p\ndiv \Delta $, where
\begin{equation} \Delta := \prod_{1\le i<j\le k}|h_j-h_i|, 
\label{6.13}
\end{equation}
since then all $k$ of the $h_i$'s are distinct modulo $p$. We now introduce an important parameter
$U$ that is used throughout the rest of the paper.  We want $U$ to be an upper bound for $\log \Delta$,
and since trivially $\Delta \le h^{k^2}$  we choose
\begin{equation} U := Ck^2\log(2h) \label{6.14} \end{equation}
and have
\begin{equation} \log \Delta \le U . \label{6.15}\end{equation}

We now prove, for
$-1/4 < \sigma \le 1$,
\begin{equation} |G_{\mathcal{H}}(s)| \ll \exp (5k U^\delta \log\log U), \quad \text{where}\ \delta
= \max (-\sigma, 0). 
\label{6.16} 
\end{equation}
We treat separately the different pieces of the product defining
$G_{\mathcal{H}}$.  First, by use of the inequality $\log(1+x) \le x$ for $x\ge 0$, we have
\[ \begin{split}  \left|\prod_{p\le U} \left(1 -
\frac{\nu_p(\mathcal{H})}{p^{1+s}}\right)\right| &\le \prod_{p\le U} \left(1 +
\frac{k}{p^{1-\delta}}\right) \\& = \exp\bigg(\sum_{p\le U}\log\Big( 1  + \frac{k}{p^{1-\delta}}\Big)\bigg) \\& \le \exp\bigg( \sum_{p\le U} \frac{k}{p^{1-\delta}}\bigg)\\ &
\le \exp\bigg( kU^\delta\sum_{p\le U}\frac{1}{p} \bigg) \\ & \ll \exp\left( kU^\delta \log\log U\right). \end{split}\]
Second, by the same estimates and the inequality $(1-x)^{-1} \le 1+3x$ for $0\le x \le 2/3$, we see that
\[ \begin{split}  \left|\prod_{p\le U} \left(1 -
\frac{1}{p^{1+s}}\right)^{-k}\right| &\le \bigg(\prod_{p\le U} \left(1 -
\frac{1}{p^{1-\delta}}\right)^{-1}\bigg)^k \\& \le \bigg(\prod_{p\le U}\left(1 +
\frac{3}{p^{1-\delta}}\right)\bigg)^k \quad \Big(\text{since} \ \frac{1}{p^{1-\delta}} \le \frac{1}{{2}^{3/4}} <\frac23\Big),\\ & \ll \exp\left( 3kU^\delta \log\log U\right). \end{split}\]
Hence, the terms in the product for $G_\mathcal{H}(s)$ with $p\le U$ are  $\ll \exp\left( 4kU^\delta \log\log U\right)$. 

For the terms $p>U$, we first consider those for which $p|\Delta$. In absolute value, they are
\[ \le \prod\limits_{\substack{p\mid \Delta \\ 
p >U}} \left(1 + \frac{k}{p^{1-\delta}}\right) \left(1 +
\frac{3}{p^{1 - \delta}} \right)^k   \ 
\le \ \exp\Big(\sum_{\substack{p\mid \Delta \\
p >U}} \frac{4k}{p^{1-\delta}}\Big) .\]
Since there are less than $(1+o(1)) \log \Delta <U$ primes with $p|\Delta$, the sum above is increased if we replace these terms with the integers between $U$ and $2U$.  Therefore the right-hand side above is
\[\le \exp\Big( 4k \sum_{U<n\le 2U} \frac{1}{n^{1-\delta}}\Big)  \le 
\exp\Big( 4k (2U)^\delta \sum_{U<n\le 2U} \frac{1}{n}\Big) 
\le \exp(4kU^\delta).  \]
Finally, if $p>U$   
\[  \left|\frac{k}{p^{1+s}}\right| \le \frac{k}{U^{1-\delta}} \le \frac{1}{2} ,\]
so that in absolute value the terms with $p>U$ and $p\ndiv \Delta$ are
\[\begin{split} &= \Bigg|\prod_{\substack{p\ndiv \Delta\\ p>U}} \left(1 -
\frac{k}{p^{1+s}}\right)
\left(1-\frac{1}{p^{1+s}}\right)^{-k}\Bigg| \\&
=\Bigg| \exp\Bigg(\sum_{\substack{p\ndiv \Delta\\ p>U}}\bigg(-\sum_{\nu=1}^\infty \frac{1}{\nu}\Big(\frac{k}{p^{1+s}}\Big)^\nu + k\sum_{\nu=1}^\infty \frac{1}{\nu}\Big(\frac{1}{p^{1+s}}\Big)^\nu\bigg)\Bigg)\Bigg|\\&
\le \exp\bigg(\sum_{p>U}\sum_{\nu=2}^\infty \frac{2}{\nu}\Big(\frac{k}{p^{1-\delta}}\Big)^\nu\bigg) \\&
\le \exp\bigg(\sum_{p>U}\sum_{\nu=2}^\infty \Big(\frac{k}{p^{1-\delta}}\Big)^\nu\bigg) \\&
\le \exp\bigg(2k^2\sum_{n>U}\frac{1}{n^{2-2\delta}}\bigg)\\&
\le \exp\Big(\frac{4k^2U^\delta}{U^{1-\delta}}\Big)\\&
\le \exp\big(2kU^\delta\big) .\end{split}\]
Thus,  the terms with $p>U$ contribute $\le \exp\big(6kU^\delta\big)$, from which we obtain \eqref{6.16}.

In conclusion,
for $ h\ll R^C$ (where $C>0$ is fixed and as large
as we wish) and for $s$ on or to the right of $\mathcal{L}$, we have
\begin{equation}
F(s) \ll (C \log(|t|+3))^k
\exp(5k U^\delta \log\log U).
\label{6.17}\end{equation}
Returning to  the integral in \eqref{6.7}, we see that the integrand vanishes as
$|t|\to \infty$, $-1/4 < \sigma \leq 1$.
By \eqref{6.9} we see that in moving  the contour from $(1)$ to the left to $\mathcal{L}$ we either pass through a simple pole at $s=0$ when $\mathcal{H}$ is admissible (so that $\gs(\mathcal{H})\neq 0$), or we pass through a regular point at $s=0$ when  $\mathcal{H}$ is not admissible.  In either case, we have by virtue of
\eqref{5.2}, \eqref{6.11}, \eqref{6.14},
\eqref{6.17}, and Lemma 1,
for any $k$ satisfying \eqref{6.1},
\begin{equation}
\begin{split}
\mathcal{T}_R(N;\mathcal{H}) &=
G_{\mathcal{H}}(0) +\frac {1}{ 2\pi
i}\mathop{\int}_{\mathcal{L}\ } F(s)\frac{R^{s}}{ {s}^{k+1}}
ds\\
&  = \gs(\mathcal{H}) +O(e^{-c\sqrt{\log R}}).
\label{6.18}
\end{split}
\end{equation}
Equation \eqref{6.2} now follows from this and \eqref{6.5}.

\begin{rema}
The  exponent $1/2$ in the restriction $k \ll (\log R)^{1/2 - \eta_0}$
 is not significant. Using Vinogradov's zero-free region for
$\zeta(s)$ we could replace $1/2$ by $3/5$.
\end{rema}

\section{First Part of the Proof of Proposition 1}
Let
\begin{equation}
\gathered \mathcal{H}=\mathcal{H}_1\cup \mathcal{H}_2, \
|\mathcal{H}_1|=k_1, \ |\mathcal{H}_2|=k_2, \ k=k_1+k_2, \\  r=
|\mathcal{H}_1\cap \mathcal{H}_2|, \  M = k_1 + k_2 + \ell_1 +
\ell_2. \endgathered \label{7.1}
\end{equation}
Thus $|\mathcal{H}|= k-r$.
We prove  Proposition~1 in the following sharper form.

\begin{Prop1prime}
Let $h\ll R^C$, where $C$ is any positive fixed constant. As $R,N\to \infty$, we have
\begin{equation}
\begin{split}
\sum_{n\le N}
\Lambda_R( n;\mathcal{H}_1,\ell_1)\Lambda_R( n;\mathcal{H}_2,\ell_2)
&=
{\ell_1 + \ell_2 \choose \ell_1} \frac{(\log R)^{r + \ell_1 +
\ell_2}}{(r + \ell_1 + \ell_2)!} \mathfrak S(\mathcal H)N \\
&\qquad + N
\sum_{j=1}^{r + \ell_1 + \ell_2}
\mathcal{D}_j(\ell_1,\ell_2,\mathcal{H}_1,\mathcal{H}_2)(\log
R)^{r+ \ell_1 + \ell_2 -j} \\
&\qquad  +O_M\left(Ne^{-c\sqrt{\log R}}\right)
+ O (R^2(3\log R)^{3k + M}), 
\label{7.2}
\end{split}
\end{equation}
where the $\mathcal{D}_j(\ell_1, \ell_2, \mathcal{H}_1,\mathcal{H}_2)$'s are
functions
independent of $R$ and $N$\marginpar{}
 which satisfy the bound
\begin{equation}
\mathcal{D}_j(\ell_1, \ell_2,\mathcal{H}_1,\mathcal{H}_2)
\ll_M  (\log U)^{C_j}
\ll_M (\log\log 10 h)^{C'_j}
\label{7.3}
\end{equation}
where $U$ is defined in \eqref{6.14} and $C_j$ and $C'_j$ are two
positive constants depending on $M$.
\end{Prop1prime}
\noindent
\textit{Proof}.  We can assume that both $\mathcal{H}_1$ and
$\mathcal{H}_2$ are non-empty since the case where one of these
sets is empty can be covered in the same way as we did in case of
$\ell = 0$ in Section 6. Thus $k\ge 2$ and
we have
\begin{equation}
\begin{split}
&\mathcal{S}_R(N; \mathcal{H}_1,\mathcal{H}_2, \ell_1, \ell_2)
:=\sum_{n=1}^N \Lambda_R( n;\mathcal{H}_1,\ell_1)\Lambda_R( n;\mathcal{H}_2,\ell_2) \\
&\hspace{0.5cm}= \frac{1}{(k_1 + \ell_1)!(k_2 + \ell_2)!}
\sum_{d,e\le R}\mu(d)\mu(e)\left(\log \frac{R}{d}\right)^{k_1 +
\ell_1} \left(\log \frac{R}{e}\right)^{k_2 + \ell_2}\sum_{\substack{1\le n\le N\\
d|P_{\mathcal{H}_1}(n)\\ e|P_{\mathcal{H}_2}(n) }}1. \label{7.4}
\end{split}\end{equation}
For the inner sum, we let $d=a_1a_{12}$, $e=a_2a_{12}$ where
$(d,e)=a_{12}$. Thus $a_1$, $a_2$, and $a_{12}$ are pairwise
relatively prime, and the divisibility conditions $
d|P_{\mathcal{H}_1}(n)$ and $e|P_{\mathcal{H}_2}(n)$ become
$a_1|P_{\mathcal{H}_1}(n)$, $a_2|P_{\mathcal{H}_2}(n)$,
$a_{12}|P_{\mathcal{H}_1}(n)$, and
$a_{12}|P_{\mathcal{H}_2}(n)$. As in Section 6, we get
$\nu_{a_1}(\mathcal{H}_1)$ solutions for $n$ modulo $a_1$, and
$\nu_{a_2}(\mathcal{H}_2)$ solutions for $n$ modulo $a_2$. If
$p|a_{12}$, then from the two divisibility conditions we have
$\nu_p(\mathcal{H}_1(p) \cap  \mathcal{H}_2(p))$ solutions  for
$n$ modulo $p$, where
\[ \mathcal{H}(p) = \{ {h'}_1,\ldots ,{h'}_{\nu_p(\mathcal{H})}:
{h'}_j\equiv h_i\in
\mathcal{H}\  \text{for some}\  i, 1\le {h'}_j \le p\}\]
Notice that $\mathcal{H}(p)= \mathcal{H}$ if $p>h$ . Alternatively,
we can avoid this definition which is necessary only for small
primes by defining
\begin{equation}
\overline \nu_p(\mathcal{H}_1 \overline \cap \thinspace \mathcal{H}_2) :=
\nu_p(\mathcal{H}_1(p) \cap  \mathcal{H}_2(p)) :=
\nu_p(\mathcal{H}_1)+ \nu_p(\mathcal{H}_2) - \nu_p(\mathcal{H})
\label{7.5}
\end{equation}
and then extend this definition to squarefree numbers by
multiplicativity.\footnote{We are making a convention here that
for $\overline \nu_p$ we take intersections modulo $p$.}
\marginpar{}
Thus we see that
\[
\sum_{\substack{1\le n\le N\\  d|P_{\mathcal{H}_1}(n)\\
e|P_{\mathcal{H}_2}(n)}} 1 = \nu_{a_1}(\mathcal{H}_1)
\nu_{a_2}(\mathcal{H}_2) \overline\nu_{a_{12}}(\mathcal{H}_1
\overline \cap\thinspace
\mathcal{H}_2)\left( \frac{N}{a_1a_2a_{12}} +O(1)\right) ,\]
\marginpar{}
and have
\begin{equation}
\begin{split}
&\mathcal{S}_R(N; \ell_1, \ell_2, \mathcal{H}_1,\mathcal{H}_2)\\
& \quad = \frac{N}{(k_1 + \ell_1)!(k_2 + \ell_2)!}
\sumprime_{\substack{ a_1a_{12}\le R\\
a_2a_{12}\le R}}\frac{\mu(a_1)\mu(a_2)\mu(a_{12})^2\nu_{a_1}
(\mathcal{H}_1)\nu_{a_2}(\mathcal{H}_2)\overline \nu_{a_{12}}
(\mathcal{H}_1 \overline \cap\thinspace
\mathcal{H}_2)}{a_1a_2a_{12}}\\
& \hskip 2in \times \left(\log
\frac{R}{a_1a_{12}}\right)^{k_1 + \ell_1} \left(\log
\frac{R}{a_2a_{12}}\right)^{k_2 + \ell_2}\\
& \hskip .35in + O\left((\log R)^M \sumprime_{\substack{ a_1a_{12}\le R\\
a_2a_{12}\le R}} \mu(a_1)^2 \mu(a_2)^2 \mu(a_{12})^2
\nu_{a_1}(\mathcal{H}_1)\nu_{a_2}(\mathcal{H}_2) \overline
\nu_{a_{12}}(\mathcal{H}_1 \overline\cap\thinspace
\mathcal{H}_2)\right)\\
&\quad = N \mathcal{T}_R(\ell_1, \ell_2;\mathcal{H}_1,\mathcal{H}_2)
+ O(R^2 (3\log R)^{3k+M}), \label{7.6}
\end{split}
\end{equation}
where $\sum'$ indicates the summands are pairwise relatively prime.
Notice that by Lemma 2, the error term was bounded by
\[ \begin{split}
&\ll (\log R)^M \sideset{}{^\flat}\sum_{q \leq R^2} \sum_{q =
a_1 a_2 a_{12}} d_k(q)\\
&= (\log R)^{M} \sideset{}{^\flat}\sum_{q\leq R^2} d_3(q) d_k(q)\\
&= (\log R)^{M} \sideset{}{^\flat}\sum_{q \leq R^2} d_{3k}(q)\\
& \ll R^2(3\log R)^{3k+M}.
\end{split} \]
By \eqref{6.6}, we have
\begin{equation}
\mathcal{T}_R(\ell_1, \ell_2; \mathcal{H}_1,\mathcal{H}_2)= \frac
{1}{ (2\pi i)^{2}} \int\limits_{(1)} \! \int\limits_{(1)}
F(s_1,s_2)\frac{R^{s_1}}{ {s_1}^{k_1+\ell_1 + 1}} \frac{R^{s_2}}{
{s_2}^{k_2+\ell_2 + 1}} \,ds_1 \,ds_2 , \label{7.7}
\end{equation}
where, by letting $s_j=\sigma_j+it_j$ and assuming $\sigma_1,
\sigma_2>0$,
\begin{equation}
\begin{split} F(s_1, s_2) &= \sumprime_{ 1\le
a_1,a_2,a_{12} <\infty} \frac{\mu(a_1)\mu(a_2) \mu(a_{12})^2
\nu_{a_1}(\mathcal{H}_1)\nu_{a_2}(\mathcal{H}_2)\overline\nu_{a_{12}}
(\mathcal{H}_1 \overline \cap \thinspace
\mathcal{H}_2)}{{a_1}^{1+s_1}{a_2}^{1+s_2}{a_{12}}^{1+s_1+s_2}}
\\&
= \prod_p\biggl(1 - \frac{\nu_p(\mathcal{H}_1)}{p^{1+s_1}} -
\frac{\nu_p(\mathcal{H}_2)}{p^{1+s_2}} + \frac{\overline
\nu_p(\mathcal{H}_1 \overline \cap \thinspace
\mathcal{H}_2)}{p^{1+s_1+s_2}}\biggr).
\end{split}
\label{7.8}
\end{equation}
Since for all $p>h$ we have $\nu_p(\mathcal{H}_1) =k_1$,
$\nu_p(\mathcal{H}_2) =k_2$, and $\nu_p(\mathcal{H}_1\cap
\mathcal{H}_2)=r$, we factor out the dominant zeta-factors and write
\begin{equation} F(s_1,s_2) =
G_{\mathcal{H}_1,\mathcal{H}_2}(s_1,s_2)
\frac{\zeta(1+s_1+s_2)^r}{\zeta(1+s_1)^{k_1} \zeta(1+s_2)^{k_2}},
\label{7.9}
\end{equation}
where by \eqref{5.1}
\begin{equation}
G_{\mathcal{H}_1,\mathcal{H}_2}(s_1,s_2) =
\prod_p\left(\frac{\left(1 - \frac{\nu_p(\mathcal{H}_1)}{p^{1+s_1}}
- \frac{\nu_p(\mathcal{H}_2)}{p^{1+s_2}} + \frac{\overline
\nu_p(\mathcal{H}_1 \overline\cap \mathcal{H}_2)}{p^{1+s_1+s_2}}
\right)\left(1-\frac{1}{p^{1+s_1+s_2}}\right)^r}
{\left(1-\frac{1}{p^{1+s_1}}\right)^{k_1}\left(1-\frac{1}
{p^{1+s_2}}\right)^{k_2}}\right) \label{7.10}
\end{equation}
is analytic and uniformly bounded for $\sigma_1, \sigma_2 > - 1/4 +
\delta$, for any fixed $\delta > 0$. Also, from \eqref{2.2},
\eqref{7.1}, and \eqref{7.5} we see immediately that
\begin{equation} G_{\mathcal{H}_1,\mathcal{H}_2}(0,0) =
\gs(\mathcal{H}) .\label{7.11}
\end{equation} 
Furthermore, the same argument leading to \eqref{6.16}  shows that
for $s_1$, $s_2$ on $\mathcal L$ or to the right of $\mathcal
L$ 
\begin{equation}
G_{\mathcal{H}_1,\mathcal{H}_2}(s_1,s_2) \ll
\exp (Ck U^{\delta_1 + \delta_2} \log\log U).
\label{7.12}
\end{equation}
with $\delta_i = -\min(\sigma_i, 0)$ and $U$ defined in \eqref{6.14}.
Thus for $s_1$ and $s_2$ on $\mathcal L$ or
to the right of $\mathcal{L}$ we have
\begin{equation}
F(s_1,s_2) \ll \exp(Ck U^{\delta_1 + \delta_2} \log\log U)
\big(\log (2+|t_1|)\log(2+ |t_2|)\big)^{2k}
\max \left(1,\frac{1}{|s_1+s_2|^r}\right) .
\label{7.13}
\end{equation}
The integrand of \eqref{7.7}\marginpar{} vanishes as either $|t_1|\to
\infty$ or $|t_2| \to \infty$, $\sigma_1, \sigma_2 \in [-c, 1]$.
We define
\begin{equation} W(s) := s\zeta(1+s) 
\label{7.14}\end{equation}
and 
\begin{equation}  D(s_1,s_2) = G_{\mathcal{H}_1,\mathcal{H}_2}(s_1,s_2)
\frac{W(s_1+s_2)^r}{W(s_1)^{k_1} W(s_2)^{k_2}},
\label{7.15}
\end{equation}
so that
\begin{equation}
\mathcal{T}_R(\ell_1, \ell_2; \mathcal{H}_1,\mathcal{H}_2)= \frac {1}{ (2\pi
i)^{2}}
\mathop{\int }_{(1)}
\mathop{\int}_{(1)}
D(s_1,s_2)\frac{R^{s_1+s_2}}{ {s_1}^{\ell_1 + 1}
{s_2}^{\ell_2 + 1}(s_1+s_2)^r}\,ds_1ds_2 .
\label{7.16}
\end{equation}
To complete the proof of Proposition 1, we need to evaluate this integral. We will also need to evaluate a similar integral in the proof of Proposition 2, where the parameters $k_1$, $k_2$,  and $r$ have several slightly different relationships with $\mathcal{H}_1$ and $\mathcal{H}_2$, and $G$ is slightly altered. Therefore we change notation to handle these situations simultaneously.
\section{Completion of the proof of Proposition 1: Evaluating an integral}
Let
\begin{equation}
 \mathcal T^*_R( a, b, d, u, v, h):
= \frac1{(2\pi i)^2}
\int\limits_{(1)} \int\limits_{(1)} \frac{D(s_1, s_2)
R^{s_1 + s_2} }{s^{u+1}_1 s^{v+1}_2 (s_1+s_2)^d}\,ds_1\, ds_2,
\label{8.1}
\end{equation}
where
\begin{equation}
D(s_1, s_2) = \frac{G(s_1, s_2) W^d(s_1 + s_2)}{W^a(s_1)W^b(s_2)}
\label{8.2}
\end{equation}
and $W$ is from \eqref{7.14}. We assume $G(s_1,s_2)$ 
is regular on $\mathcal L$ and to the right of $\mathcal L$ and satisfies the bound 
\begin{equation}
G(s_1,s_2) \ll_M
\exp (CM U^{\delta_1 + \delta_2} \log\log U), \quad \text{where} \ \ U= CM^2\log(2h). 
\label{8.3}
\end{equation} 
\begin{lemma} Suppose that
\begin{equation}  0\le a,b,d,u,v\le M, \quad a+u \ge 1, \quad b+v\ge 1, \quad d\le \min(a,b) ,\label{8.4} \end{equation}
where $M$ is a large constant and our estimates may depend on $M$. 
Let $h\ll R^C$, with $C$ any positive fixed constant. Then we have,  as $R \to \infty$,
\begin{equation}
\begin{split}\mathcal T^*_R( a, b, d, &u, v, h):
=
{u+v \choose u} \frac{(\log R)^{u+v+d}}{(u+v+d)!} G(0,0) \\
& + 
\sum_{j=1}^{u+v+d}
\mathcal{D}_j(a,b,d,u,v,h)(\log
R)^{u+v+d -j}  +O_M\left(e^{-c\sqrt{\log R}}\right), 
\label{8.5}
\end{split}
\end{equation}
where the $\mathcal{D}_j(a,b,d,u,v,h)$'s are
functions
independent of $R$ \marginpar{}
 which satisfy the bound
\begin{equation}
\mathcal{D}_j(a,b,d,u,v,h)
\ll_M  (\log U)^{C_j}
\ll_M (\log\log 10 h)^{C'_j}
\label{8.6}
\end{equation}
for some positive constants $C_j$, $C'_j$ depending on $M$.
\end{lemma}

\emph{Proof.}
As in Section 7, we see
the integrand in \eqref{8.1} vanishes as $|t_1|
\to \infty$ or $|t_2| \to \infty$.
We first shift the contour $(1)$ for the integral over $s_1$ to $\mathcal L$,
passing a pole at $s_1 = 0$, and obtain
\begin{equation}
\mathcal T^*_R = \frac{1}{2\pi i} \int\limits_{(1)}\!
\mathop{\mathrm{Res}}_{s_1 = 0} \biggl(\frac{D(s_1, s_2) R^{s_1 +
s_2}}{s^{u+1}_1 s^{v+1}_2 (s_1+s_2)^d} \biggr) \,ds_2  + \frac{1}{(2\pi i)^2} \int\limits_{(1)}\!\!
\int\limits_{\mathcal L}\! \frac{D(s_1, s_2) R^{s_1 +
s_2}}{s^{u+1}_1 s^{v+1}_2 (s_1+s_2)^d} \,ds_1 \,ds_2. \label{8.7}
\end{equation} 
In the first term, we move the contour over $s_2$ along $(1)$  to $\mathcal{L}$, and pass a pole at $s_2=0$. For the second term, after interchanging the order of integration, 
we move the contour $(1)$ to the left to
$\mathcal{L}$ passing poles at $s_2= -s_1$ and 
$s_2=0$.  We thus obtain
\begin{equation}
\begin{split}
&T^*_R = \mathop{\mathrm{Res}}_{s_2=0}\mathop{\mathrm{Res}}_{s_1=0}\frac{D(s_1, s_2) R^{s_1 + s_2}}{s^{u+1}_1 s^{v+1}_2
(s_1+s_2)^d} + \frac {1}{ 2\pi i}\mathop{\int}_{\mathcal{L} }
\mathop{\mathrm{Res}}_{s_1 = 0}
\bigg(\frac{D(s_1, s_2)
R^{s_1 + s_2} }{s^{u+1}_1 s^{v+1}_2 (s_1+s_2)^d}
\bigg)\, ds_2\\& \ + \frac {1}{ 2\pi i}\mathop{\int}_{\mathcal{L} }
\mathop{\mathrm{Res}}_{s_2 = 0}
\bigg(\frac{D(s_1, s_2)
R^{s_1 + s_2} }{s^{u+1}_1 s^{v+1}_2 (s_1+s_2)^d}
\bigg)\, ds_1  +\frac {1}{ 2\pi
i}\mathop{\int}_{\mathcal{L} }
\mathop{\mathrm{Res}}_{s_2=-s_1}
\bigg(\frac{D(s_1, s_2)
R^{s_1 + s_2} }{s^{u+1}_1 s^{v+1}_2 (s_1+s_2)^d}
\bigg)\, ds_1\\
& \  + \frac {1}{ (2\pi
i)^{2}}\mathop{\int }_{\mathcal{L} }\mathop{\int}_{\mathcal{L}\
} \frac{D(s_1, s_2)
R^{s_1 + s_2} }{s^{u+1}_1 s^{v+1}_2 (s_1+s_2)^d} \, ds_1 \, ds_2  := I_0 + I_1 + I_2+ I_3 + I_4.
\end{split}
\label{8.8}
\end{equation}
We will see that the residue $I_0$ provides the main term and some of the lower order terms, the integral $I_3$ provides the remaining lower order terms, and the integrals $I_1$, $I_2$, and $I_4$ are error terms.

We consider first $I_0$.  At $s_1=0$ there is a pole of order  $\le u+1$, and therefore \footnote{If $G(0,0)=0$ then the order of the pole is $u$ or less, but the formula we use to compute the residue is still valid. In this situation one or more of the initial terms will have the value zero.}  by Leibniz's rule we have
\[\mathop{\mathrm{Res}}_{s_1=0}\frac{D(s_1, s_2) R^{s_1 }}{s^{u+1}_1 
(s_1+s_2)^d} = \frac1{u!} 
\sum^u_{i=0} {u\choose i} (\log R)^{u-i}
\frac{\partial^i}{\partial s^i_1} \left( \frac{D(s_1, s_2)}{(s_1
+ s_2)^d}\right) \Bigg|_{s_1 = 0}  \]
and 
 \[ \begin{split} \frac{\partial^i}{\partial s^i_1} &\left( \frac{D(s_1, s_2)}{(s_1
+ s_2)^d}\right) \Bigg|_{s_1 = 0}= (-1)^i \frac{D(0, s_2)
d(d+1) \cdots (d + i - 1)}{s^{d+i}_2}
\\& \quad 
+ \sum^i_{j = 1} {i\choose j} \frac{\partial^j}{\partial s^j_1} D(s_1, s_2)
\Bigg|_{s_1 = 0}  (-1)^{i-j} \frac{d(d+1) \cdots (d + i - j -
1)}{s^{d+i - j}_2} ,\end{split} \]
where in case of $i = j$ (including the case when $i = j = 0$ and $d \geq  0$ arbitrary)
 the empty product in the numerator is $1$. We conclude that
\begin{equation} \mathop{\mathrm{Res}}_{s_1=0}\frac{D(s_1, s_2) R^{s_1 }}{s^{u+1}_1 
(s_1+s_2)^d} = \sum_{i=0}^u\sum_{j=0}^i \frac{a(i,j)(\log R)^{u-i} }{s^{d+i - j}_2}\frac{\partial^j}{\partial s^j_1} D(s_1, s_2)
\Bigg|_{s_1 = 0}  \label{8.9}\end{equation}
with $a(i,j)$ given explicitly in the previous equations. 
To complete the evaluation of $I_0$,  we  see that the $(i,j)$th term contributes to $I_0$ a pole at  $s_2=0$ of order
$v + 1 + d + i - j$ (or less), and therefore by Leibniz's formula
\[ \begin{split}&\mathop{\mathrm{Res}}_{s_2=0}\frac{R^{s_2 }}{s^{v+1+d+i-j }_2}\frac{\partial^j}{\partial s^j_1} D(s_1, s_2)
\Bigg|_{s_1 = 0} \\& =\frac{1}{(v+d+i-j)!} \sum_{m=0}^{v+d+i-j}{v+d+i-j \choose m} (\log
R)^{v+d+i-j-m}  \frac{\partial^m }{\partial s^m_2}
\frac{\partial^j}{\partial s^j_1} D(s_1, s_2) \Bigg|_{\substack{s_1 =0\\ s_2
= 0}}    .\end{split}\]
This completes the evaluation of $I_0$, and we conclude
\begin{equation}
I_0 =  \sum^u_{i=0} \sum_{j=0}^i\sum_{m=0}^{v+d+i-j} b(i,j,m)\bigg(\frac{\partial^m }{\partial s^m_2}
\frac{\partial^j}{\partial s^j_1} D(s_1, s_2)  \Bigg|_{\substack{s_1 =0\\ s_2
= 0} }\ \bigg) (\log
R)^{u+v+d-j-m}, 
\label{8.10} \end{equation}
where
\begin{equation} b(i,j,m)=
(-1)^{i-j}{u\choose i} 
{i\choose j}   {v+d+i-j \choose m}\frac{d(d+1)\cdots (d+i - j -
1)}{u!(v+d+i-j)!} .
\label{8.11}
\end{equation}

The main term is of order $(\log R)^{u+v+d}$ and occurs when $j=m=0$. Therefore, it is given by
\[
 G(0,0) (\log R)^{u+v+d}\left(\frac{1}{u!}
\sum^u_{i=0} {u\choose i}
(-1)^i \frac{d(d+1)\cdots (d + i - 1)}{(v+d+i)!} \right).\]
It is not hard to prove that
\begin{equation} \frac{1}{u!}
\sum^u_{i=0} {u\choose i}
(-1)^i \frac{d(d+1)\cdots (d + i - 1)}{(v+d+i)!} = {u+v\choose u}\frac{1}{(u+v+d)!} ,
\label{8.12}\end{equation}
from which we conclude that the main term is
\begin{equation} G(0,0){u+v\choose u}\frac{1}{(u+v+d)!}(\log R)^{d+u+v}.\label{8.13} \end{equation}
Motohashi found the following approach which avoids proving \eqref{8.12} directly and can be used to simplify some of the previous analysis. Granville also made a similar observation. The residue we are computing is equal to
\[ \frac1{(2\pi i)^2}
\int\limits_{\Gamma_2} \int\limits_{\Gamma_1} \frac{D(s_1, s_2)
R^{s_1 + s_2} }{s^{u+1}_1 s^{v+1}_2 (s_1+s_2)^d}\, ds_1\, ds_2,\]
where $\Gamma_1$ and $\Gamma_2$ are the circles $|s_1|=\rho$ and $|s_2|=2\rho$, respectively, with a small $\rho>0$. Writing $s_1=s$ and $s_2=sw$, this is equal to 
\[ \frac1{(2\pi i)^2}\int\limits_{\Gamma_3} 
\int\limits_{\Gamma_1} \frac{D(s, sw)
R^{s(w+1)}}{s^{u+v+d+1} w^{v+1}(w+1)^d}\,  ds\, dw,\]
with $\Gamma_3$ the circle $|w|=2$. The main term is obtained from the constant term $G(0,0)$ in the Taylor expansion of $D(s,sw)$ and, therefore, equals
\[ G(0,0) \frac{(\log R)^{u+v+d}}{(u+v+d)!} \ \frac{1}{2\pi i} \int\limits_{\Gamma_3} \frac{(w+1)^{u+v}
}{w^{v+1}}\, dw = G(0,0) \frac{(\log R)^{u+v+d}}{(u+v+d)!} {u+v\choose v},\]
by the binomial expansion. 

To complete the analysis of $I_0$, we only need to show that the partial derivatives of $D(s_1,s_2)$ at $(0,0)$ satisfy the bounds given in the lemma. For this, we use
Cauchy's estimate for derivatives which takes the form, for $z = \sigma + it$ and $z_0 = \sigma_0 + it_0$,
\begin{equation}
|f^{(j)}(z_0)| \leq  \max_{|z - z_0| = \eta} |f(z)| \frac{j!}{ \eta^{j}},
\label{8.14}
\end{equation}
if $f(z)$ is analytic for $|z - z_0| \leq \eta$.
In the application below we will choose 
\begin{equation}
\eta = \frac{1}{C \log U \log T}, \quad \text{where} \ T = |s_1|+|s_2| + 3.
\label{8.15}
\end{equation}
We see that if $z_0$ is on $\mathcal L$ or to the right of $\mathcal
L$ then the whole circle $|z - z_0-1| = \eta$ will remain in
the region \eqref{5.3} and the estimates \eqref{5.4} hold in this circle. (We remind the reader that the generic constants $c,C$
take different values at different appearances.) Thus, we have
for $s_1, s_2$ on $\mathcal L$ or
to the right of $\mathcal L$,
\begin{equation}
\aligned
&\frac{\partial^m}{\partial
s^m_2} \frac{\partial^j}{\partial s^j_1} D(s_1, s_2) \\
&\le j! m! (C \log U)^{j+m} (\log T_1)^j (\log T_2)^m
\max_{|s_1^* - s_1| \leq \eta, |s_2^* - s_2| \leq \eta}
|D(s_1^*, s_2^*)| \\
& \ll_M \frac{\exp(CMU^{\delta_1 + \delta_2} \log\log U)
(\log T_1)^{M} (\log T_2)^{M} }
{\max (c, |s_1|)^a \max(c, |s_2|)^b}
\max(1, |s_1 + s_2|)^d ,
\endaligned
\label{8.16}
\end{equation}
which, if $\max (|s_1|, |s_2|) \leq C$, reduces to
\begin{equation}
 \frac{\partial^m}{\partial
s^m_2}\frac{\partial^j}{\partial s^j_1}  D(s_1, s_2) \ll_M  \exp\left(CM U^{\delta_1 +
\delta_2} \log\log U\right).
\label{8.17}
\end{equation}
In particular, we have
\begin{equation}
 \frac{\partial^m}{\partial
s^m_2}\frac{\partial^j}{\partial s^j_1}  D(s_1, s_2)\Bigg|_{\substack{s_1=0\\s_2=0}} \ll_M (\log U)^{C(M)}. 
\label{8.18} \end{equation}
We conclude from \eqref{8.10}, \eqref{8.13}, and \eqref{8.18} that $I_0$ provides the main term and some of the secondary terms in Lemma 3 which satisfy the stated bound. 

We now consider $I_1$. By \eqref{8.9}, \eqref{8.16}, and Lemma 1, we have

\begin{equation}\begin{split}
I_1&\ll_M (\log R)^u\int\limits_{\mathcal L} \frac{ \max(1, |s_2|)^d (\log(|t_2| +
3))^{M} e^{CMU^{\delta_2} \log\log U} R^{-\delta_2}}{|s_2|^{v+1+b+d}}
|ds_2| \\
&\ll_M \int\limits_{\mathcal L} \frac{e^{CMU^{\delta_2}\log\log U}
(\log(|t_2| + 3))^{M} R^{-\delta_2}}{|s_2|^{v+1+b} }|ds_2| \ll_M  e^{-c\sqrt{\log R}}.
\end{split}
\label{8.19}
\end{equation}
The same bound holds for $I_2$ since it is with relabeling equal to $I_1$. Further, $I_4$ also satisfies this bound by \eqref{7.13} (with relabeling) and Lemma 1. 

Finally, we  examine $I_3$, which only occurs if $d\ge 1$. Because
\begin{equation}
\begin{split}\mathop{\mathrm{Res}}_{s_2=-s_1}\bigg(
\frac{D(s_1, s_2)
R^{s_1 + s_2} }{s^{u+1}_1 s^{v+1}_2 (s_1+s_2)^d}
\bigg) & = \lim_{s_2\to
-s_1}\frac{1}{(d-1)!}\frac{\partial^{d-1}}{{\partial
s_2}^{d-1}}\left(\frac{
D(s_1,s_2)   R^{s_1+s_2}}{{s_1}^{u+1}
{s_2}^{v+1} }\right) \\
&= \frac1{(d-1)!}
\sum^{d-1}_{i=0} \mathcal{B}_i(s_1)(\log
R)^{d-1-i},
\end{split}
\label{8.20}
\end{equation}
where
\begin{equation}
\mathcal{B}_i(s_1) =
{d-1\choose i} \sum^i_{j = 0} {i\choose j}
\frac{\partial^{i-j}}{\partial s^{i-j}_2} D(s_1,
s_2)\Bigg|_{s_2 = -s_1}  
\frac{(-1)^{j} (v+1) \cdots (v+j)}{(-1)^{j + v + 1} s^{u+v+j+2}_1},
\label{8.21}
\end{equation}
we have
\begin{equation}
I_3 = \frac{1}{(d-1)!}
 \sum^{d-1}_{i=0}
\mathcal{C}_i(\log R)^{d -1-i},
\label{8.22}
\end{equation}
where
\begin{equation}
\mathcal{C}_i =\frac {1}{ 2\pi
i}\int\limits_{\mathcal{L}
}\mathcal{B}_i(s_1)\, ds_1,\ \
0\le i\le d - 1 .
\label{8.23}
\end{equation}
It remains to estimate the $\mathcal{C}_i$'s, which are independent
of $R$ but depend on $h$.  

By \eqref{8.21} we see
that the functions $\mathcal B_i(s_1)$ tend to zero as $|t_1|
\to \infty$, $-c \leq \sigma \leq 1$, and  further by  \eqref{8.16}
\begin{equation}
\mathcal{B}_i
\ll_M
\sum^i_{j=0}
 (\log T_1)^{M} \exp\left(CMU^{2\delta_1} \log\log U\right)
\frac{1} {|t_1|^{ u + v + j + 2 + a + b}}.
\label{8.24}
\end{equation}
Therefore, we may shift the contour $\mathcal{L}$ back
to the imaginary axis with a semicircle of radius $1/\log U$
centered and to the left of $s_1=0$.
The contribution to $\mathcal{C}_i$  from the
integral along the imaginary axis is
\begin{equation}
\ll_M  (\log U)^{u+v+i+a+b+1} \exp( CM\log\log  U) \ll_M (\log U)^{C'(M)}.
\label{8.25}
\end{equation}
This expression also bounds the
contribution to $\mathcal{C}_i$  from  the
semicircle contour 
and thus completes the evaluation of $I_3$. Combining our results, we obtain Lemma 3.
\section{Proof of Proposition 2}
We introduce some standard notation associated with \eqref{1.2} and \eqref{1.3}. Let
\begin{equation} \theta(x;q,a) := \sum_{\substack{p\le x \\  p\equiv a (\text{mod}\, q)}}\log p  = [ (a,q)=1]\frac{x}{\phi(q)} + E(x;q,a),\label{9.1} \end{equation}
where $[S]$ is 1 if the statement $S$ is true and is $0$ if $S$ is false. Next, we define
\begin{equation} E'(x,q) := \max_{\substack{a\\  (a,q)=1}} |E(x;q,a)|,  \quad  E^*(x,q) = \max_{y\le x}E'(y,q). \label{9.2} \end{equation}
In this paper we only need level of distribution results for $E'$, but usually these results are stated in the stronger form for $E^*$. Thus, for some $1/2\le \vartheta\le 1$, we assume, given any $A>0$ and $\epsilon>0$, that
\begin{equation} \sum_{q\le x^{\vartheta - \epsilon}} E^*(x,q) \ll_{A,\epsilon} \frac{x}{(\log x)^A}.  \label{9.3} \end{equation}
This is known to hold with $\vartheta =1/2$. 

We prove the following stronger version of Proposition~2.  Let
\begin{equation}
C_R(\ell_1, \ell_2, \mathcal{H}_1,\mathcal{H}_2,h_0) = \left\{
\begin{array}{ll}
      {1} &
        \mbox{if $h_0\not \in \mathcal{H}$,}\\
 \frac{(\ell_1 + \ell_2 + 1) \log R}{(\ell_1 + 1)(r + \ell_1 + \ell_2 +1)}
 &
        \mbox{if $h_0 \in \mathcal{H}_1$ and $h_0\not \in
 \mathcal{H}_2$,}\\
 \frac{(\ell_1 + \ell_2 + 2)(\ell_1 + \ell_2 + 1) \log R}
{(\ell_1 + 1)(\ell_2 + 1)(r + \ell_1 + \ell_2 + 1)} &
        \mbox{if $h_0 \in \mathcal{H}_1\cap \mathcal{H}_2$.}\\
\end{array}
\right.\label{9.4}
\end{equation}
By relabeling the variables we obtain the corresponding form
if $h_0 \in \mathcal H_2$ and $h_0 \not\in \mathcal H_1$. We continue to use the notation \eqref{7.1}.
\begin{Prop2prime}
Suppose  $h\ll R$.  Given any positive $A$,
there is a $B = B(A,M)$ such that for $R \ll_{M,A}
N^{\frac{1}{4}}/(\log N)^B$ and $R,N\to \infty$,
\begin{equation}
\begin{split}
\sum_{n=1}^N&\Lambda_R(n;\mathcal{H}_1,\ell_1)
\Lambda_R(n;\mathcal{H}_2, \ell_2)\theta(n+h_0)
\\
&= \frac{C_R(\ell_1, \ell_2, \mathcal{H}_1,\mathcal{H}_2,h_0) }
{(r + \ell_1 + \ell_2)!}
{\ell_1 + \ell_2 \choose \ell_1}
\gs({\mathcal{H}}^0) N (\log R)^{r + \ell_1 + \ell_2}\\
&\ \ + N
\sum_{j=1}^r\mathcal{D}_j(\ell_1, \ell_2,
\mathcal{H}_1,\mathcal{H}_2,h_0)(\log
R)^{r + \ell_1 + \ell_2 - j} +O_{M,A}\left(\frac{N}{(\log N)^A}\right),
\label{9.5}
\end{split}
\end{equation}
where the $\mathcal{D}_j(\ell_1, \ell_2,
\mathcal{H}_1,\mathcal{H}_2,h_0)$'s are
 functions independent of $R$ and $N$ which satisfy the bound
\begin{equation}
\mathcal{D}_j(\mathcal{H}_1,\mathcal{H}_2,h_0)
\ll_M (\log U)^{C_j} \ll_M(\log\log 10 h)^{C'_j}
\label{9.6}
\end{equation}
for some positive constants $C_j$, $C'_j$ depending on $M$.
Assuming that conjecture \eqref{9.3}  holds, then \eqref{9.5} holds for $R\ll_{M}
N^{\frac{\vartheta}{2} -\epsilon}$ and $h\le R^\epsilon$,  for any given
$\epsilon>0$.
\end{Prop2prime}

\noindent
\textit{Proof}.
We will assume that both $\mathcal{H}_1$ and $\mathcal{H}_2$ are
non-empty so that $k_1\ge 1$ and $k_2 \ge 1$. The proof in the case when one of
these sets is empty is much easier and may be obtained by
an argument analogous to that of Section~6.
We have
\begin{equation}
\aligned
&\widetilde{\mathcal{S}}_R(N;\mathcal{H}_1,\mathcal{H}_2, \ell_1,
\ell_2, h_0) 
:=\sum_{n=1}^N
\Lambda_R(n;\mathcal{H}_1, \ell_1)\Lambda_R(n;\mathcal{H}_2, \ell_2)
\theta(n+h_0)\\
&= \frac{1}{(k_1 + \ell_1)!(k_2 + \ell_2)!}\sum_{d,e\le
R}\mu(d)\mu(e)\left(\log \frac{R}{d}\right)^{k_1 + \ell_1} \left(\log
\frac{R}{e}\right)^{k_2 + \ell_2}\sum_{\substack{1\le n\le N\\
d|P_{\mathcal{H}_1}(n)\\ e|P_{\mathcal{H}_2}(n)
}}\theta(n+h_0).
\endaligned
\label{9.7}
\end{equation}
To treat the inner sum above,  let $d=a_1a_{12}$ and $e=a_2a_{12}$, where
$(d,e)=a_{12}$, so that $a_1$, $a_2$, and $a_{12}$ are pairwise
relatively prime. As in Section 7, the $n$ for which
$d|P_{\mathcal{H}_1}(n)$ and  $e|P_{\mathcal{H}_2}(n)$ cover\marginpar{}
certain residue classes modulo $[d,e]$. If  $n\equiv b \,
(\text{mod} \, a_1a_2a_{12})$ is such a residue class, then
letting
$m=n+h_0 \equiv b+h_0(\text{mod} \, a_1a_2a_{12})$,  we see that this
residue class contributes to the inner sum
\begin{equation}
\begin{split}
&\sum_{\substack{1+h_0 \leq  m\le N+h_0\\ m\equiv b
+h_0\, (\text{mod} \, a_1a_2a_{12}) }}\theta(m)
= \theta(N+h_0;a_1a_2a_{12},b+h_0) -
\theta(h_0;a_1a_2a_{12},b+h_0)\\& \qquad =
[(b+h_0,a_1a_2a_{12})=1]\frac{N}{\phi(a_1a_2a_{12})} +E(N;
a_1a_2a_{12},b+h_0)+O(h\log N).\end{split}
\label{9.8}
\end{equation}
We need to determine the number of these residue classes where
$(b+h_0,a_1a_2a_{12})=1$ so that the main term is non-zero.
If $p|a_1$, then $b\equiv - h_j\ (\text{\rm mod}\, p)$ for some $h_j \in
\mathcal{H}_1$, and therefore
$b+h_0\equiv h_0 - h_j\ (\text{\rm mod}\,
p)$.  Thus, if $h_0$ is distinct modulo $p$ from all the $h_j\in
\mathcal{H}_1$, then all $\nu_p(\mathcal{H}_1)$ residue classes
satisfy the relatively prime condition, while otherwise
$h_0\equiv h_j(\text{\rm mod}\, p)$ for some $h_j\in \mathcal{H}_1$
leaving $\nu_p(\mathcal{H}_1)-1$ residue classes with a non-zero
main term. We introduce the notation $\nu_p^*({\mathcal{H}_1}^0)$
 for this number in either case,
where we define for a set $\mathcal{G}$ and integer $h_0$
\begin{equation}
\nu_p^*(\mathcal{G}) = {\nu_p}(\mathcal{G}^0)-1, \label{9.9}
\end{equation}
where
\begin{equation}  \mathcal{G}^0 = \mathcal{G} \cup  \{h_0\}.
\label{9.10}\end{equation}  
 We extend this definition to
$\nu_d^*({\mathcal{H}_1}^0)$ for squarefree numbers $d$ by
multiplicativity. The function $\nu_d^*$ is familiar in sieve
theory, see \cite{HR}. A more algebraic discussion of  ${\nu_d}^*$ may also be found in \cite{GGPY,GMPY}. We define $\overline \nu^*_d \big((\mathcal H_1 \overline \cap
\mathcal H_2)^0\big)$  as in \eqref{7.5}.

Next, the divisibility conditions  $a_2|P_{\mathcal{H}_2}(n)$,
$a_{12}|P_{\mathcal{H}_1}(n)$, and
$a_{12}|P_{\mathcal{H}_2}(n)$ are handled as in Section 7
together with the above considerations. Since
$E(n;q,a)\ll (\log N)$ if $(a,q)>1$ and $q \leq N$, we conclude that
\begin{equation}
\begin{split} &\sum_{\substack{1\le n\le N\\
d|P_{\mathcal{H}_1}(n)\\ e|P_{\mathcal{H}_2}(n)}}
\theta(n+h_0) = \nu_{a_1}^*({\mathcal{H}_1}^0)
\nu_{a_2}^*({\mathcal{H}_2}^0)\overline\nu_{a_{12}}^*
\left((\mathcal{H}_1 \overline \cap
\mathcal{H}_2)^0\right)\frac{N}{\phi(a_1a_2a_{12})} \\
& \quad + O\left( d_k(a_1a_2a_{12})\left(\max_{\substack{b \\ (b,
a_1a_2a_{12})=1}}\big|E(N ;a_1a_2a_{12}, b)\big|\ + h(\log
N)\right) \right).\end{split}\label{9.11}
\end{equation}\marginpar{}
Substituting this into \eqref{9.7} we obtain
\begin{equation}
\begin{split}
&\widetilde{\mathcal{S}}_R(N;\mathcal{H}_1,\mathcal{H}_2, \ell_1, \ell_2,
h_0)\\
&
\quad  = \frac{N}{(k_1 + \ell_1)!(k_2 + \ell_2)!}
\sumprime_{\substack{ a_1a_{12}\le
R\\ a_2a_{12}\le
R}}\frac{\mu(a_1)\mu(a_2)\mu(a_{12})^2\nu_{a_1}^*({\mathcal{H}_1}^0)
\nu_{a_2}^*({\mathcal{H}_2}^0)\overline
\nu_{a_{12}}^*\left((\mathcal{H}_1 \overline \cap
\mathcal{H}_2)^0\right)}{\phi(a_1a_2a_{12})} \\&
\hskip 2in
\times \left(\log \frac{R}{a_1a_{12}}\right)^{k_1 + \ell_1} \left(\log
\frac{R}{a_2a_{12}}\right)^{k_2 + \ell_2}\\
&  + O\left((\log R)^M
\sumprime_{\substack{ a_1a_{12}\le R\\ a_2a_{12}\le R}}
d_k(a_1a_2a_{12})E'(N,a_1a_2a_{12})\right)  +
O(hR^2 (3\log N)^{M+3k + 1})\\
&
= N\widetilde{ \mathcal{T}}_R(\mathcal{H}_1,\mathcal{H}_2,
\ell_1, \ell_2, h_0)
+O \left((\log R)^M\mathcal{E}_k(N)\right)+ O(hR^2 (3\log
N)^{M+3k+1}), \label{9.12}
\end{split}\end{equation}
where the last error term was obtained using Lemma 2.
To estimate the first error term we use
Lemma~2, \eqref{1.3},  and the trivial estimate
$E'(N,q)\leq  (2N/q)\log
N$ for $q\le N$, and \eqref{9.3} to find,
uniformly for $k \leq \sqrt{(\log N)/18}$, that
\begin{equation}
\begin{split}
|\mathcal{E}_k(N)| &\leq  \sideset{}{^\flat}\sum_{q\le R^2} d_k(q)
\max_{\substack{b \\ (b, q)=1}}
\big|E(N ;  q,
b)\big|\sum_{q=a_1a_2a_{12}}1
\\
&
=  \sideset{}{^\flat}\sum_{q\le R^2} d_k(q) d_3(q)
E'(N, q)\\
&
\leq \sqrt{\sideset{}{^\flat}\sum_{q\le
R^2} \frac{d_{3k}(q)^2}{q}}
\sqrt{\sum_{q\le R^2}q
(E'(N, q))^2}\\
&
\leq \sqrt{(\log N)^{9k^2}}\sqrt{2N\log N}
\sqrt{\sum_{q\le
R^2} E'(N,q)} \\
&
\ll  N(\log N)^{(9k^2 + 1 - A)/2},
\end{split}
\label{9.13}
\end{equation}
provided $R^2 \ll N^{\frac{1}{2}}/(\log N)^B$. On relabeling, we conclude that given any positive integers $A$ and
$M$ there is a positive constant $B=B(A,M)$ so that for  $R\ll
\frac{N^{\frac{1}{4}}}{(\log N)^B} $ and $h\le R$,
\begin{equation}
\widetilde{\mathcal{S}}_R(N;\mathcal{H}_1,\mathcal{H}_2, \ell_1, \ell_2, h_0)=
N\widetilde{\mathcal{T}}_R(\mathcal{H}_1,\mathcal{H}_2,\ell_1, \ell_2, h_0)+
O_M\left(\frac{N}{(\log N)^A}\right).
\label{9.14}
\end{equation}
Using \eqref{9.3} with any $\vartheta>1/2$, we see that \eqref{9.14} holds for the longer range
$R\ll_M N^{\frac{\vartheta}{2}-\epsilon}$, $h \ll N^\epsilon$.

Returning to the main term in \eqref{9.12}, we have by \eqref{6.6} that
\begin{equation}
\widetilde{\mathcal{T}}_R(\mathcal{H}_1,\mathcal{H}_2, \ell_1, \ell_2, h_0)= \frac {1}{
(2\pi i)^{2}}\mathop{\int }_{(1)}
\mathop{\int}_{(1)}
F(s_1,s_2)
\frac{R^{s_1}}{ {s_1}^{k_1 + \ell_1 +1}}
\frac{R^{s_2}}{{s_2}^{k_2 + \ell_2 +1}}ds_1ds_2 ,\label{9.15}
\end{equation}
where, by letting $s_j=\sigma_j+it_j$ and assuming $\sigma_1,
\sigma_2>0$,
\begin{equation}
\begin{split} F(s_1, s_2) &= \sumprime_{ 1\le
a_1,a_2,a_{12} <\infty}
\frac{\mu(a_1)\mu(a_2) \mu(a_{12})^2
\nu_{a_1}^*({\mathcal{H}_1}^0)\nu_{a_2}^*
({\mathcal{H}_2}^0)\overline\nu_{a_{12}}^*
\left(({\mathcal{H}_1 \overline \cap
\mathcal{H}_2})^0\right)}
{\phi(a_1){a_1}^{s_1}\phi(a_2){a_2}^{s_2}\phi(a_{12}){a_{12}}^{s_1+s_2}}
\\&
=\prod_p\bigg(1 -
\frac{\nu_p^*({\mathcal{H}_1}^0)}{(p-1)p^{s_1}} -
\frac{\nu_p^*({\mathcal{H}_2}^0)}{(p-1)p^{s_2}}
+ \frac{\overline \nu_p^*\left(({\mathcal{H}_1 \overline\cap
\mathcal{H}_2})^0\right)}{(p-1)p^{s_1+s_2}}\bigg).\end{split}\label{9.16}
\end{equation}
We now consider three cases.

\medskip
\emph{Case 1.} Suppose $h_0 \not \in  \mathcal{H}$. Then we have, for
$p>h$,
\[ \nu_p^*({\mathcal{H}_1}^0) =k_1, \quad
\nu_p^*({\mathcal{H}_2}^0) =k_2, \quad
\overline\nu_p^*\left(({\mathcal{H}_1 \overline\cap
\mathcal{H}_2})^0\right)=
r.\]
Therefore in this case we define
$G_{\mathcal{H}_1,\mathcal{H}_2}(s_1,s_2)$  by
\begin{equation} F(s_1,s_2) =
G_{\mathcal{H}_1,\mathcal{H}_2}(s_1,s_2)
\frac{\zeta(1+s_1+s_2)^r}{\zeta(1+s_1)^{k_1}\zeta(1+s_2)^{k_2}}.
\label{9.17} \end{equation}

\medskip
\emph{Case 2.} Suppose $h_0 \in  \mathcal{H}_1$  but $h_0 \not
\in \mathcal{H}_2$. (By relabeling this also covers the case
where $h_0 \in  \mathcal{H}_2$  and $h_0 \not \in
\mathcal{H}_1$.) Then for $p>h$
\[ \nu_p^*({\mathcal{H}_1}^0) =k_1-1, \quad
\nu_p^*({\mathcal{H}_2}^0) =k_2, \quad
\overline\nu_p^*\left(({\mathcal{H}_1\overline
\cap \mathcal{H}_2})^0\right)=
r.\]
Therefore, we define $G_{\mathcal{H}_1,\mathcal{H}_2}(s_1,s_2)$  by
\begin{equation} F(s_1,s_2) =
G_{\mathcal{H}_1,\mathcal{H}_2}(s_1,s_2)
\frac{\zeta(1+s_1+s_2)^r}{\zeta(1+s_1)^{k_1-1}\zeta(1+s_2)^{k_2}}.
\label{9.18} \end{equation}

\medskip
\emph{Case 3.} Suppose $h_0  \in  \mathcal{H}_1\cap
\mathcal{H}_2$. Then for $p>h$ 
\[ \nu_p^*({\mathcal{H}_1}^0) =k_1-1, \quad
\nu_p^*({\mathcal{H}_2}^0) =k_2-1, \quad
\overline\nu_p^*\left(({\mathcal{H}_1 \overline\cap
\mathcal{H}_2})^0\right)=
r-1.\]
Thus, we define
$G_{\mathcal{H}_1,\mathcal{H}_2}(s_1,s_2)$  by
\begin{equation} F(s_1,s_2) =
G_{\mathcal{H}_1,\mathcal{H}_2}(s_1,s_2)
\frac{\zeta(1+s_1+s_2)^{r-1}}{\zeta(1+s_1)^{k_1-1}
\zeta(1+s_2)^{k_2-1}}.\label{9.19}
\end{equation}
In each case,  $G$ is analytic and uniformly bounded for $\sigma_1, \sigma_2 >
-c$, with any $c<1/4$.

We now show that in all three cases
\begin{equation} G_{\mathcal{H}_1,\mathcal{H}_2}(0,0) =
\gs(\mathcal{H}^0) .\label{9.20}
\end{equation}
Notice that in the second two cases we have $\mathcal{H}^0 =
\mathcal{H}$.   By \eqref{5.1}, \eqref{7.5}, \eqref{9.9}, and \eqref{9.16}, we find in all three cases
\begin{equation}
\begin{split}
&G_{\mathcal{H}_1,\mathcal{H}_2}(0,0) \\
&\quad =
\prod_p\left( 1 -
\frac{\nu_p({\mathcal{H}_1}^0)+\nu_p({\mathcal{H}_2}^0)
- \overline \nu_p((\mathcal{H}_1 \overline\cap \mathcal{H}_2)^0)-1}{p-1}
\right)\left(1-\frac{1}{p}\right)^{-a(\mathcal{H}_1,\mathcal{H}_2,h_0)}\\&
\quad = \prod_p\left( 1 -
\frac{\nu_p({\mathcal{H}}^0)-1}{p-1}\right)\left(
1-\frac{1}{p}\right)^{-a(\mathcal{H}_1,\mathcal{H}_2,h_0)},\end{split}
\label{9.21}
\end{equation}
where in Case 1 \
$a(\mathcal{H}_1,\mathcal{H}_2,h_0)=k_1+k_2-r =
k-r$, in Case 2 \
$a(\mathcal{H}_1,\mathcal{H}_2,h_0)=(k_1-1)+k_2-r = k-r-1$, and
in Case 3 \
$a(\mathcal{H}_1,\mathcal{H}_2,h_0)=(k_1-1)+(k_2-1)-(r-1) =
k-r-1$.  Hence, in Case~1 we have
\begin{equation}
 \begin{split}G_{\mathcal{H}_1,\mathcal{H}_2}(0,0) &
= \prod_p\left(
\frac{p-\nu_p({\mathcal{H}}^0)}{p-1}\right)\left(
1-\frac{1}{p}\right)^{-(k-r)}\\&=\prod_p\left( 1 -
\frac{\nu_p({\mathcal{H}}^0)}{p}\right)\left(
1-\frac{1}{p}\right)^{-(k-r+1)} \\&= \gs({\mathcal{H}}^0),
\end{split}
\label{9.22}
\end{equation}
while in Cases 2 and 3 we have
\begin{equation}
\begin{split}G_{\mathcal{H}_1,\mathcal{H}_2}(0,0) &
= \prod_p\left(  \frac{p-\nu_p(\mathcal{H})}{p-1}\right)\left(
1-\frac{1}{p}\right)^{-(k-r-1)}\\&=\prod_p\left( 1 -
\frac{\nu_p(\mathcal{H})}{p}\right)\left(
1-\frac{1}{p}\right)^{-(k-r)} \\&= \gs(\mathcal{H}) \quad (
=\gs({\mathcal{H}}^0) ) .\end{split}
\label{9.23}
\end{equation}

We are now ready to evaluate
$\mathcal{T}_R(\mathcal{H}_1,\mathcal{H}_2,\ell_1, \ell_2, h_0)$. There are
two differences between  the functions $F$ and $G$ that appear in
\eqref{9.16}--\eqref{9.19} and the earlier \eqref{7.8}--\eqref{7.10}.
The first difference is that a factor
of $p$ in the denominator of the Euler product in \eqref{7.8}
has  been replaced by $p-1$, which only effects the value of
constants in calculations.
The second difference is the
relationship between $k_1$, $k_2$, and $r$, which effects the
residue calculations  of the main terms.
However, the analysis of
lower order terms and the error analysis is essentially unchanged
 and,
therefore, we only need to examine the main terms. We use Lemma 3 here to
cover all of the cases.
Taking into account \eqref{9.17}--\eqref{9.19} we have in Case 1 that $a=k_1,b=k_2, d=r, u=\ell_1, v=\ell_2$; in Case 2 that $a=k_1-1,b=k_2, d=r, u=\ell_1+1, v=\ell_2$; and in Case 3 that $a=k_1-1,b=k_2-1, d=r-1, u=\ell_1+1, v=\ell_2+1$. 
By \eqref{9.22} and \eqref{9.23}, the proof of Theorems~2 and 2$'$ is thus complete.

\section{Proof of Theorem 3}
For convenience, we agree in our notation below that we
consider every set of size $k$ with a multiplicity $k!$
according to all permutations of the elements $h_i \in \mathcal
H$, unless mentioned otherwise. While unconventional, this will clarify some of the calculations. 

To prove Theorem 3 we consider in place of \eqref{3.5}
\begin{equation}\begin{split}
&\mathcal{S}_R(N,k,\ell,h,\nu) \\&
:= \frac{1}{N h^{2k+1}}\sum_{n=N+1}^{2N}
\bigg( \sum_{1\le h_0\le h} \theta (n+h_0) - \nu \log 3N \bigg)\Bigg(
\sum_{\substack{\mathcal{H} \subset \{1,2,\ldots, h\} \\
|\mathcal{H}|=k}}
\Lambda_{R}(n;\mathcal{H}, \ell )\Bigg)^2 \\&
= \widetilde M_R(N,k,\ell,h)  - \nu \frac{\log 3N}{ h} M_R(N,k,\ell,h),
\label{10.1}
\end{split}\end{equation}
where 
\begin{equation}
M_R(N,k, \ell,h) = \frac{1}{Nh^{2k}} \sum^{2N}_{n = N+1}
\Bigg(
\sum_{\substack{\mathcal{H} \subset \{1,2,\ldots, h\} \\
|\mathcal{H}|=k}}
\Lambda_{R}(n;\mathcal{H}, \ell )\Bigg)^2 
\label{10.2}
\end{equation}
and
\begin{equation}
\widetilde M_R(N,k,\ell,h) = \frac{1}{Nh^{2k+1}} \sum^{2N}_{n = N+1}
\bigg(\sum_{1\le h_0\le h} \theta(n+h_0) \bigg) \Bigg(
\sum_{\substack{\mathcal{H} \subset \{1,2,\ldots, h\} \\
|\mathcal{H}|=k}}
\Lambda_{R}(n;\mathcal{H}, \ell )\Bigg)^2 .
\label{10.3}
\end{equation}
To evaluate $M_R$ and $\widetilde M_R$ we multiply out the sum and apply Propositions 1 and 2.  We need to group the pairs of
sets $\mathcal H_1$ and $\mathcal H_2$ according to size of the intersection $r= |\mathcal H_1 \cap \mathcal H_2|$, and thus $|\mathcal H| = |\mathcal H_1 \cup \mathcal
H_2| = 2k - r$.
Let us choose now a set $\mathcal H$ and here exceptionally we
disregard the permutation of the elements in $\mathcal H$.
(However for $\mathcal H_1$ and $\mathcal H_2$ we take into account
all permutations.)
Given the set $\mathcal H$ of size $2k - r$, we can
choose $\mathcal H_1$ in ${2k - r\choose k}$ ways. Afterwards, we
can choose the intersection set in ${k\choose r}$ ways. Finally,
we can arrange the elements both in $\mathcal H_1$ and $\mathcal
H_2$ in $k!$ ways. This gives
\begin{equation}
{2k - r\choose k} {k \choose r} (k!)^2  = (2k - r)! {k \choose r}^2
 r!
\label{10.4}
\end{equation}
choices for $\mathcal H_1$ and $\mathcal H_2$, taking into
account the permutation of the elements in $\mathcal H_1$  and $\mathcal H_2$.
If we consider in the summation every union set $\mathcal H$
of size $j$ just once,
independently of the arrangement of the elements, then
Gallagher's theorem \eqref{3.7} may be formulated as
\begin{equation}\sideset{}{^*}\sum_{\substack{\mathcal{H} \subset \{1,2,\ldots, h\} \\
|\mathcal{H}|=j}}
\
\mathfrak S(\mathcal H)\ \sim \  \frac{h^j}{j!}
,
\label{10.5}
\end{equation}
where $\sum^*$ indicates every set  is counted just once.
Applying this, we obtain on letting
\begin{equation} x = \frac{\log R}{h} ,\label{10.6} \end{equation} 
and using Proposition 1
\begin{equation}
\aligned
M_R(N,k,\ell,h)  &\sim \frac1{Nh^{2k}} \sum^k_{r=0} (2k - r)! {k
\choose r}^2 r!  {2\ell \choose \ell}
\frac{(\log R)^{2\ell + r}}{(r + 2\ell)!}  N\sum_{|\mathcal
H| = 2k - r}{}^{^{\hspace*{-3.5mm}\scriptstyle *}}
\hspace*{3mm} \mathfrak S(\mathcal H)\\
&\sim {2\ell \choose \ell} (\log R)^{2\ell} \sum^k_{r=0}
{k\choose r}^2  \frac{x^r} {(r + 1) \dots (r + 2\ell)}.
\endaligned
\label{10.7}
\end{equation}
By Proposition 2 and \eqref{10.5} we have
\begin{equation}
\widetilde M_R (N, k, \ell,h)  \sim
\frac1{{N}h^{2k+1}} \sum^k_{r=0} (2k - r)!
{k\choose r}^2 \, r! \, Z_r,
\label{10.8}
\end{equation}
where abbreviating $a = \frac{2\ell + 1}{\ell + 1} = {2\ell + 1
\choose \ell + 1} {2\ell \choose \ell}^{-1} =  \frac{1}{2}{2\ell + 2\choose
\ell + 1} {2\ell \choose \ell}^{-1} $, we have
\begin{equation}\begin{split}
Z_r : &= {2\ell \choose \ell} \frac{(\log R)^{2\ell + r}}{(r +
2\ell)!} \Bigg\{ r \sideset{}{^*}\sum_{|\mathcal H| = 2k -
r}{}2a \frac{\log
R}{r + 2\ell + 1} \mathfrak S(\mathcal H)N\\
&+ (2k - 2r) \sideset{}{^*}\sum_{|\mathcal H| = 2k -
r}
 a \frac{\log
R}{r + 2\ell + 1} \mathfrak S(\mathcal H) N +
\sideset{}{^*}\sum_{|\mathcal H| = 2k - r} 
\sum^h_{h_0 = 1 \atop h_0 \notin {\mathcal H}}
\mathfrak S(\mathcal H^0)N\Bigg\}\\
& \sim N {2\ell \choose \ell} \frac{(\log R)^{2\ell + r}}{(r +
2\ell)!}
\left\{ \frac{h^{2k - r}}{(2k - r)!}  \frac{2a k \log R}{r
+ 2\ell + 1} + \frac{2k - r + 1}{(2k - r + 1)!} h^{2k - r + 1}
\right\},
\label{10.9}
\end{split}
\end{equation}
where in the last sum we took into account which element of
$\mathcal H^0$ is $h_0$, which can be chosen in $2k - r + 1$
ways.
Thus we obtain
\begin{equation}
\widetilde M_R(N,k,\ell,h) \sim {2\ell \choose \ell} (\log R)^{2\ell}
\sum^k_{r = 0} {k\choose r}^2 \frac{x^r}{(r + 1) \dots (r
+ 2\ell)} \left(\frac{2 a k}{r + 2\ell + 1} x + 1 \right).
\label{10.10}
\end{equation}
We conclude, on introducing the parameters
\begin{equation}  \varphi =\frac{1}{\ell+1}, \ (\text{so}\  a= 2- \varphi), \quad  \Theta = \frac{\log R}{\log 3N}, \ (\text{so} \  R= (3N)^\Theta), \label{10.11} \end{equation}
that
\begin{equation}
S_R(N, k, \ell, h,\nu)  \sim {2\ell \choose \ell} (\log R)^{2\ell} P_{k, \ell,\nu}(x),
\label{10.12}
\end{equation}
where
\begin{equation}
P_{k,\ell,\nu}(x) = \sum^k_{r = 0} {k\choose r}^2 \frac{x^r}{(r+1) \cdots (r
+ 2\ell)} \left(1 + x\left(\frac{4(1 - \frac{\varphi}{2})k}{r +
2\ell + 1} - \frac{\nu}{\Theta}\right)\right).
\label{10.13}
\end{equation}
Let
\begin{equation}  h=  \lambda \log 3N,  \quad \text{so that}\ \  x = \frac{\Theta}{\lambda} . \label{10.14}\end{equation}
The analysis of when $S>0$ now depends on the polynomial $P_{k,\ell,\nu}(x)$.  We examine this polynomial as $k, \ell \to \infty$ in such a way that $\ell = o(k)$.  In the first place, the size of the terms of the polynomial are determined by the factor 
\[ g(r) ={k\choose r}^2 x^r ,\]
and since $g(r)>g(r-1)$  is equivalent to 
\[ r< \frac{k+1}{1 +\frac{1}{\sqrt{x}}}\]
we should expect that the polynomial is controlled by terms with $r$ close to  $k/(z+1)$, where
\begin{equation} z = \frac{1}{\sqrt{x}}. \label{10.15} \end{equation}
Consider now the sign of each term.  For small $x$, the terms in the polynomial are positive, but they become negative when 
\[ 1 + x\left(\frac{4(1 - \frac{\varphi}{2})k}{r +
2\ell + 1} - \frac{\nu}{\Theta}\right) <0 .\]
When $r= k/(z+1)$ and letting  $k,\ell \to \infty$, $\ell =o(k)$,  we have heuristically
\[ \begin{split} 1 + x\left(\frac{4(1 - \frac{\varphi}{2})k}{r +
2\ell + 1} - \frac{\nu}{\Theta}\right) &\approx 1 + \frac{1}{z^2}\left( \frac{4k}{\frac{k}{z+1}} -\frac{\nu}{\Theta}\right)\\& = \frac{1}{z^2}\big( (z+2)^2 - \frac{\nu}{\Theta}\big).\end{split} \]
Therefore, the terms will be positive for $r$ in this range if $ z< \sqrt{\frac{\nu}{\Theta}} -2$, which is equivalent to $
\lambda < (\sqrt{ \nu} - 2\sqrt{\Theta})^2$. Since we can take $\Theta$ as close to $\vartheta/2$ as we wish, this will imply Theorem 3. 
To make this argument precise, we choose $r_0$ slightly smaller than where $g(r)$ is maximal, and prove that all the negative terms together contribute less then the single term $r_0$, which will be positive for $z$ and thus $\lambda$ close to the values above.

For the proof, we may assume $\nu \geq 2$ and $1/2\le \vartheta_0 \le 1$ are fixed, with
$\vartheta_0 < 1$ in case of $\nu = 2$. (The case $\nu=1$ is covered by Theorem 2, and 
the case $\nu = 2$, $\vartheta_0 = 1$, $E_2 = 0$ is covered by \eqref{1.11}  proved in Section 3.)
First, we choose $\varepsilon_0$ as a sufficiently small fixed positive number.
We will choose $\ell$ sufficiently large, depending on $\nu$,
$\vartheta_0$, $\varepsilon_0$, and set
\begin{equation}
k = (\ell+1)^2 = \varphi^{-2}, \quad
\ell > \ell_0(\nu, \vartheta_0, \varepsilon_0),\ \ \text{so that}\
\varphi < \varphi_0(\nu, \vartheta_0, \varepsilon_0).
\label{10.16}
\end{equation}
Furthermore, we choose
\begin{equation}
\Theta = \frac{\log R}{\log 3N} = \frac{\vartheta_0 (1 -
\varphi)}{2},
\label{10.17}
\end{equation}
and (because of our assumptions on $\nu$) we can define
\begin{equation}
z_0 :=  \sqrt{2\nu / \vartheta_0} - 2>0.
\label{10.18}
\end{equation}
Thus, we see that
\begin{equation}  1 + \frac{1}{{z_0}^2}\left( \frac{4k}{\frac{k}{z_0+1}} -\frac{2\nu}{\vartheta_0}\right) = \frac{1}{{z_0}^2}\big( (z_0+2)^2 - \frac{2\nu}{\vartheta_0}\big)=0,\label{10.19} \end{equation}
 Let us choose now
\begin{equation}
r_0 = \left[ \frac{k + 1}{z_0 + 1} \right], \quad
r_1 = r_0 + \varphi k = r_0 + \ell+1,
\label{10.20}
\end{equation}
and put  
\begin{equation} z = z_0(1 + \varepsilon_0) .\label{10.21}\end{equation} 
The linear factor in each term of $P_{k, \ell, \nu}(x)$ is, for
 $r_0 \leq r \leq r_1$,
\begin{equation}\begin{split}
1 + x\left(\frac{4(1 - \frac{\varphi}{2})k}{r + 2\ell + 1} -
\frac{\nu}{\Theta} \right)
&= 1 + \frac{1}{z^2_0(1 + \varepsilon_0)^2} \left( \frac{4k(1 +
O(\varphi))}{\frac{k}{z_0 + 1} + O(k\varphi)} -
\frac{2\nu}{\vartheta_0 (1 - \varphi)} \right) \\& =  1 + \frac{-z^2_0 + O(\sqrt{\nu}\varphi) +
O (\nu\varphi)}{z^2_0 (1 + \varepsilon_0)^2}\\&
> c(\nu, \vartheta_0) \varepsilon_0\ \  \text{ \ if \ }
\varphi < \varphi_0(\nu, \vartheta, \varepsilon_0),
\label{10.22}
\end{split}
\end{equation}
where $c(\nu, \vartheta_0)>0$ is a constant.
Letting 
\begin{equation}
f(r) := {k\choose r}^2 \frac{x^r}{(r + 1) \dots (r + 2\ell)} ,
\label{10.23}
\end{equation}
we have, for any $r_2 > r_1$, 
\begin{align}
\frac{f(r_2)}{f(r_0)}
&< \prod_{r_0 < r \leq r_2} \left( \frac{k
+ 1 - r}{r} \cdot \frac{1}{z} \right)^2 < \prod_{r_0 < r \leq
r_1} \left( \frac{k + 1 - r}{r} \cdot \frac{1}{z} \right)^2
\label{10.24}\\
&< \left(\left( \frac{k + 1}{r_0 + 1} - 1\right) \frac{1}{z}
\right)^{2\ell} \leq \left( \frac{z_0 + 1 - 1}{z_0(1 + \varepsilon_0)} \right)^{2\ell} < e^{-\varepsilon_0 \ell}.
\nonumber
\end{align}
Thus, the total contribution in absolute value of the negative terms of $P_{k,\ell,\nu}(x)$ will be, for
sufficiently large $\ell$, at most
\begin{equation}
k\left(1 + \frac{4(k + \nu)}{z^2} \right) e^{-\varepsilon_0 \ell}
f(r_0) < e^{-\varepsilon_0 \ell / 2} f(r_0),
\label{10.25}
\end{equation}
while that of the single term $r_0$ will be by \eqref{10.22}
at least
\begin{equation}
c(\nu, \vartheta_0) \varepsilon_0 f(r_0) > e^{-\varepsilon_0
\ell/2} f(r_0) \ \text{ if } \ell > \ell_0(\nu, \vartheta_0,
\varepsilon_0).
\label{10.26}
\end{equation}
This shows that $P_{k,\ell,\nu}(x) > 0$. Hence, we must have at least $\nu+1$
primes in some interval 
\begin{equation}
[n + 1, n + h] = [n + 1, n + \lambda \log 3N], \quad n \in [N +
1, 2N],
\label{10.27}
\end{equation}
where 
\begin{equation}
\lambda  = \Theta z^2 < \frac{\vartheta_0}{2} z^2_0 (1
+ \varepsilon_0)^2 = (1 + \varepsilon_0)^2 \left(\sqrt{\nu} -
\sqrt{2\vartheta_0} \right)^2.
\label{10.28}
\end{equation}
Since $\varepsilon_0$ can be chosen arbitrarily small, this proves
Theorem~3.

\end{document}